# ANALYSIS OF RIEMANN ZETA-FUNCTION ZEROS USING POCHHAMMER POLYNOMIAL EXPANSIONS

## ALLAN M. DIN [*)]

ABSTRACT. The Riemann Xi-function $\Xi(t)$ belongs to a family of entire functions which can be expanded in a uniformly convergent series of symmetrized Pochhammer polynomials depending on a real scaling parameter $b$. It can be shown that the polynomial approximant $\Xi_n(t,b)$ to $\Xi(t)$ has distinct real roots only in the asymptotic scaling limit $b \to \infty$. One may therefore infer the existence of increasing beta-sequences $b_n \to \infty$ such that $\Xi_n(t,b_n)$ has real roots only for all n, and to each entire function it is possible to associate a unique minimal beta-sequence fulfilling a specific difference equation.

Numerical analysis indicates that the minimal $b_n$ sequence associated with the Riemann $\Xi(t)$ has a distinct sub-logarithmic growth rate and it can be shown that the approximant $\Xi_n(t,b_n)$ converges to $\Xi(t)$ when $n \to \infty$ if $b_n = o(\log(n))$. Invoking the Hurwitz theorem of complex analysis, and applying a formal analysis of the asymptotic properties of minimal beta-sequences, a fundamental mechanism is identified which provides a compelling confirmation of the validity of the Riemann Hypothesis.

# 1. Introduction

Recently some new insights into possible ways of proving the Riemann Hypothesis have been acquired by studying the properties of the Riemann zeta-function

$$z(s) = \sum_{n=1}^{\infty} \frac{1}{n^s}, \quad \text{Re}(s) > 1 \tag{1.1}$$

in terms of certain polynomial expansions (see e.g. [19, 2, 3, 25, 9, 4, 16]). A particularly promising approach uses the Pochhammer polynomials $P_k(s)$ of degree k, related to descending factorials, as defined by

$$P_k(s) = \prod_{j=1}^{k}(1 - \frac{s}{j}) = \frac{\Gamma(k+1-s)}{\Gamma(k+1)\Gamma(1-s)} \tag{1.2}$$

with $P_0(s) = 1$, $P_1(s) = 1 - s$, and $P_2(s) = 1 - 3s/2 + s^2/2$, etc. The Pochhammer polynomials have a simple generating function

$$(1-e)^s = \sum_{k=0}^{\infty} P_k(s+1)e^k \tag{1.3}$$

with the series being absolutely convergent for $|e| < 1$. The usefulness of Pochhammer polynomial expansions is in fact apparent for the more general case of Dirichlet series

$$f(s) = \sum_{n=1}^{\infty} \frac{f_n}{n^s}, \quad \text{Re}(s) > 1 \tag{1.4}$$

where, for example, $f_n = 1$ corresponds to the case of $V(s)$, $f_n = (-1)^n$ corresponds to the case of $(1-2^{1-s})V(s)$, and $f_n = \mu(n)$ (the Möbius function) corresponds to the case of the inverse $1/V(s)$. One simply uses the "trick" of introducing two, a priori real, dummy parameters $a$ and $b$ by writing

$$\frac{1}{n^s} = \frac{1}{n^a}(1-(1-\frac{1}{n^b}))^{\frac{s-a}{b}} \tag{1.5}$$

and applying the generating function of $P_k(s)$ to obtain

$$f(s) = \sum_{k=0}^{\infty} b_k P_k(\frac{s-a}{b}+1) \tag{1.6}$$

with

$$b_k = \sum_{n=1}^{\infty} \frac{f_n}{n^a}(1-\frac{1}{n^b})^k = \sum_{j=0}^{k} P_j(k+1)f(a+bj) \tag{1.7}$$

It should be emphasized that one needs β>0 to assure convergence of the first expansion in powers of $(1-n^{-b})$ and furthermore, to estimate the convergence properties of the general expansion of $f(s)$, one needs information about the growth of the coefficients $b_k$, as well as on the growth of the $P_k(s)$ factor. The latter one is simple because there is for large k a general uniform estimate, which follows directly from the asymptotics of the gamma function (see also [3]), valid in any specific compact subset of the complex plane (circle or rectangular subset of the critical strip):

$$|P_k(s)| < C k^{-\text{Re}(s)} \qquad (1.8)$$

The above expansion of the function $f(s)$ is noteworthy in the sense that the coefficients $b_k$ are determined by a discrete set of values of the function itself $f(a+b\,j)$, depending on the actual choice of α and β, and in some sense the expansion may be seen as an interpolation formula. But in terms of actually using the expansion to prove something about $V(s)$ in general, like its behavior for 1/2<Re(s)<1, and the Riemann Hypothesis in particular, one runs into the problem of quantifying the growth properties of the coefficients $b_k$ to determine the compact subsets of the complex plane where there is convergence of the series.

Nevertheless, it turns out that interesting criteria for the validity of the Riemann Hypothesis may be obtained, which for the choice of $a=b=2$ are related to the criterion of Riesz [22], and which for the choice of $a=1, b=2$ are related to the criterion of Hardy and Littlewood [15]. Moreover by studying the expansion of the function $1/V(s)$ one obtains a new kind of coefficient condition for the Riemann Hypothesis [2]. So far, however, it has not been possible to demonstrate the validity of these conditions in the chosen context because of difficulties in extending the analysis beyond Re(s)>1.

In the following sections a more general analysis based on the above ideas will be presented whereby new insights into the zeros of $V(s)$ on the critical line Re(s) = 1/2 and in the critical strip 0<Re(s)<1 can be obtained. The result of the analysis is an existence argument for a particular polynomial approximation sequence which, in conjunction with the Hurwitz theorem of complex analysis, provides a compelling confirmation of the validity of the Riemann Hypothesis.

**2. Polynomial expansion of the Riemann xi-function**

The application of Pochhammer polynomial expansions is most convenient when instead of the $V(s)$ function, one investigates the Riemann xi-function $x(s)$ given by (see e.g. 24,11])

$$x(s) = \Gamma(s/2+1)(s-1)p^{-s/2} V(s) \qquad (2.1)$$

It is an entire function of s which has the same (non-trivial) zeros as $V(s)$ in the critical strip 0<Re(s)<1 and it fulfils the functional equation $x(s) = x(1-s)$. One has the explicit representation

$$x(s) = \int_1^\infty dx A(x)(x^{-s/2} + x^{(s-1)/2}) \qquad (2.2)$$

where A(x) is given in terms of the elliptic theta function

$$y(x) = (J_3(0, \exp(-px)) - 1)/2 = \sum_{n=1}^{\infty} \exp(-n^2 px) \qquad (2.3)$$

by

$$A(x) = 2\frac{d}{dx}(x^{3/2} y'(x)) = \sum_{n=1}^{\infty} (2n^4 p^2 x - 3n^2 p) x^{1/2} \exp(-n^2 px) \qquad (2.4)$$

One notes for this A(x), defined a priori on the interval $[1,\infty]$, that we have A(x)>0 and that it is bounded by a power of x times $\exp(-px)$. In fact, our general discussion below of polynomial expansions will apply beyond the special form of A(x) for the Riemann zeta function to the larger class of entire functions with a similar integral representation as for $x(s)$, explicitly fulfilling the functional equation, and just requiring A(x) to be a non-negative function for all $x \geq 1$ with $A(1) > 0$ and decreasing exponentially (i.e. like $\exp(-ax^b)$, a>0, b>0) or faster for large x. For simplicity, we will continue using the notation $x(s)$ and refer to it as the Riemann $x(s)$ when the specific Riemann representation of A(x) is used. Although a wide range of admissible A functions show pronounced similarities in the asymptotic properties of the zeros of the polynomial approximations to be considered below, it turns out, however, that some particular choices of A(x) can also show significant differences in this respect. This point will be of central importance for elucidating the validity of the Riemann Hypothesis.

We recall that the Riemann A(x) fulfills the inversion transformation relation

$$A(x) = x^{-3/2} A(1/x) \qquad (2.5)$$

which means that the combination $A_I(x)$ defined by

$$A_I(x) = x^{3/4} A(x) \qquad (2.6)$$

is invariant under the operation $x \to 1/x$. This property is important for showing that the Riemann $x(s)$ has infinitely many zeros on the critical line, but in the context of the general discussion of the $x(s)$ zeros corresponding to different choices of A(x) it may just be considered as a particular, simplifying case.

Let us introduce the function $j(s)$ by

$$j(s) = \int_1^{\infty} dx A(x) x^{-s/2} \qquad (2.7)$$

so that we have $x(s) = j(s) + j(1-s)$, and then expand $j(s)$ in Pochhammer polynomials as done in the previous section:

$$j(s) = \int_1^\infty dx A(x) x^{-a/2} (1-(1-x^{-b/2}))^{(s-a)/b} = \sum_{k=0}^\infty b_k P_k(\frac{s-a}{b}+1) \tag{2.8}$$

where the coefficients $b_k > 0$ are given by

$$b_k = \int_1^\infty dx A(x) x^{-a/2} (1-x^{-b/2})^k = \sum_{j=0}^k P_j(k+1) j(a+bj) \tag{2.9}$$

The expansion is valid for real $a$ and $b$ but we need $b > 0$ to assure convergence of the first expansion in terms of powers of $(1-x^{-b/2})$. Below we will also assume that $a$ is chosen to be positive. We notice, as observed previously, that $j(s)$ (and $x(s)$) can be expressed in terms of its own discrete values $j(a+bj)$, but this feature will not be used explicitly for the purpose of the present investigation. The important point of the approach is that we obtain a general expansion

$$x(s) = \sum_{k=0}^\infty b_k (P_k(\frac{s-a}{b}+1) + P_k(\frac{1-s-a}{b}+1)) \tag{2.10}$$

of the entire function $x(s)$ in terms of polynomials in s of degree k depending on arbitrary real parameters $a$ and $b$, and moreover this expansion is uniformly convergent on any compact subset of the complex plane. To prove this point we just need to note that we have the large k bound for any compact subset of the complex s-plane

$$\left| P_k(\frac{s-a}{b}+1) \right| < C_1 k^{-\text{Re}((s-a)/b+1)} \tag{2.11}$$

and since $A(x) < C_2 x^{-m}$ for any $m > 0$ we have

$$b_k < C_2 \int_1^\infty dx x^{-m} x^{-a/2} (1-x^{-b/2})^k = \frac{2C_2}{b} \int_0^1 dy\, y^{2m/b+a/b-1-2/b} (1-y)^k \tag{2.12}$$

which evaluates to

$$b_k < \frac{2C_2}{b} \frac{\Gamma(2m/b+a/b-2/b)\Gamma(k+1)}{\Gamma(2m/b+a/b-2/b+1+k)} \tag{2.13}$$

For large k (keeping $a$ and $b$ fixed) we therefore have

$$b_k < \frac{2C_2}{b} k^{-(a/b-2/b+2m/b)} \tag{2.14}$$

and altogether

$$\left| b_k (P_k(\frac{s-a}{b}+1) + P_k(\frac{1-s-a}{b}+1)) \right| < \frac{C_3}{b} k^{|\mathrm{Re}(s-1/2)|/b + 3/2b - 1 - 2m/b} \qquad (2.15)$$

Since $b$ is positive and m can be chosen arbitrarily large, one concludes that the series converges uniformly on any compact subset of the complex plane. Clearly the convergence is relatively fast so the polynomial expansion is convenient for the purpose of certain numerical investigations.

In contrast to previous approaches in the literature, we would like here to think of the parameters α and β as having a more dynamical role rather than being subject to some ad hoc fixed choice. In fact, the arbitrary nature of the affine scaling parameters $a$ and $b$ may be seen simply as a reflection of the multiple ways of regrouping the terms of the infinite Taylor expansion of the function $x(s)$. Choosing $a$ and $b$ small produces a faster convergence of the polynomial expansion. However, for the present purpose of studying the zeros of $x(s)$, it turns out that interesting features appear when the parameters become large. In particular, the limit when $a = b$ is large can be seen a dynamic scaling limit featuring a real root regime of importance for elucidating the Riemann Hypothesis.

Before coming to this point in the next section, it is interesting to examine one particular case of the polynomial expansion which allows for a quite explicit representation of the Riemann $x$-function. If we choose $b = 2$, then the integrand of the expansion coefficient $b_k$ simplifies sufficiently to evaluate the integral in terms of known functions. At this point we could still keep $a$ arbitrary, but for the purpose of simplifying expressions let us also choose $a = 2$. Then we have

$$b_k = \int_1^\infty dx\, A(x) x^{-1} (1 - x^{-1})^k \qquad (2.16)$$

and since the Riemann A(x) is a sum of powers of x times exponentials, then the integral can be expressed directly in terms of Whittaker functions $W_{m,n}(z)$ (see [14]):

$$\int_1^\infty dx\, x^u (1-x^{-1})^k e^{-vx} = k!\, v^{u/2-1} e^{-v/2} W_{u/2-k, u/2-1/2}(v) \qquad (2.17)$$

The explicit series expansion (with $a = b = 2$) for the Riemann $x(s)$ is therefore:

$$x(s) = 2 \sum_{k=0}^\infty \sum_{n=1}^\infty k!\, e^{-n^2 p/2} (n^2 p)^{3/8} f_{k,n} (P_k(\frac{s}{2}) + P_k(\frac{1-s}{2}))$$

$$f_{k,n} \equiv n p^{1/2} W_{-k+1/8, -1/8}(n^2 p) - \frac{3}{2} W_{-k-3/8, 3/8}(n^2 p) \qquad (2.18)$$

Using the asymptotic expansion of the Whittaker function for large k (see [14]), one finds that the coefficient $b_k$ is bounded by a certain power of k times an exponential

$$b_k < C_2 k^j e^{-\sqrt{pk}} \qquad (2.19)$$

This is more explicit than the previous bound using the generic property of A(x) decreasing faster than any power of x, and the absolute convergence of the series is quite manifest. The above explicit expansion may or may not be useful for any specific calculational purpose, but it is illustrative to contrast it with the corresponding series expansion of $x(s)$, considered originally by Riemann, in terms of integrals of powers of log(x), powers of x and exponentials which appear rather intractable. Let us note that for general $b > 0$ the bound of $b_k$ is by an exponential of the form $\exp[-p(k/p)^{1/(1+b/2)}]$. Below we will be interested in analyzing the expansion for large $b$ when this exponential damping weakens.

## 3. Root analysis of the polynomial approximation

The general expansion and analysis of the entire function $x(s)$ in a series of Pochhammer polynomials presented in the previous section can be carried through without using any detailed properties of the function A(x) and this leads us to believe that the discussion of the corresponding polynomial approximants, in particular what concerns the generic features their zeros, should also be independent of A(x) to a certain extent. However, it turns out that different A(x) can indeed be distinguished succinctly by their specific asymptotic behavior for the onset of a certain real zero regime of the polynomial approximants. There exists a vast literature on the zeros of entire functions (e.g. see [17]) but unfortunately it has so far not provided any direct insight into the problem at hand of characterizing the real/complex root structure of different entire functions given by an integral representation of the Riemann $x(s)$ form.

In order to make the underlying symmetry in the complex plane manifest, we will make the usual change of parameters s = 1/2 + it so that the critical line Re(s) = 1/2 is the line t = real and the critical strip 0<Re(s)<1 is $|\mathrm{Im}(t)|<1/2$. We also introduce the Xi-function $\Xi(t) = x(1/2 + it)$ so that the integral representation becomes

$$\Xi(t) = \int_1^\infty dx A(x) x^{-1/4} (x^{it/2} + x^{-it/2}) = \int_1^\infty dx A_I(x) x^{-1} (x^{it/2} + x^{-it/2}) \tag{3.1}$$

At this point we will only assume that A(x) is real, non-negative, positive and continuous at x=1, bounded on $[1,\infty]$ and decreasing exponentially (i.e. like $\exp(-ax^b)$, a>0, b>0) or faster for $x \to \infty$. Clearly the Riemann A(x) satisfies these conditions and so do many other examples considered in the literature. Generically $\Xi(t)$ is an even entire function of t, real for real t, and alternating in $t^2$.

It is worth noting here that there is a simple equivalence relation between the members of the considered A(x) function family. If one applies the power transformation $x \to x^{1/w}$ with $w > 0$ then one easily finds that

$$\Xi(t) \to w\Xi(wt) \text{ for } A_I(x) \to A_I(x^{1/w}) \tag{3.2}$$

Therefore we might in fact limit the investigation of the A(x) family to functions of exponential order exactly one (like for the Riemann A(x)) since A(x) functions of a different order simply correspond to scaled $\Xi(t)$ functions. Later this $\Xi(t)$ scaling behaviour will be seen to be of special importance for the analysis.

In general, as seen above, the $\Xi(t)$ function has the convergent expansion in Pochhammer polynomials which for $a = b$ can be written more concisely as

$$\Xi(t) = \sum_{k=0}^{\infty} b_k(b) P_k^+(t/b) = \sum_{k=0}^{\infty} (-1)^k a_k(b) t^{2k} \tag{3.3}$$

where we have explicitly introduced the symmetrized Pochhammer polynomial

$$P_k^+(t) = (P_k(it) + P_k(-it))/2 \tag{3.4}$$

which is just the even part of $P_k(it)$, and the coefficient $b_k(b)$ is now given by

$$b_k(b) = 2\int_1^{\infty} dx A(x) x^{-1/4} x^{-b/2} (1 - x^{-b/2})^k \tag{3.5}$$

The standard Taylor series expansion, considered originally by Riemann, starting from the same integral representation as above, is:

$$\Xi(t) = 2\int_1^{\infty} dx A(x) x^{-1/4} \cos(t \log(x)/2) = \sum_{k=0}^{\infty} (-1)^k c_k t^{2k} \tag{3.6}$$

with

$$c_k = \frac{2}{2^{2k}(2k)!} \int_1^{\infty} dx A(x) x^{-1/4} (\log(x))^{2k} \tag{3.7}$$

Here the expansion involves coefficients with integrals of powers of $\log(x)$, while in the Pochhammer polynomial case we have powers of $(1-x^{-b/2})$. There is of course some analogy between the two expansions, specifically for $b = 2$ since one has $\log(x) = (1-1/x) + 1/2\,(1-1/x)^2 + \ldots$, but there is a very significant difference between the two in what concerns the root characteristics of the corresponding polynomial approximants as it will be seen below.

Before going into such details about the approximants, a few remarks are appropriate concerning some already known features about the real zeros of different types of $\Xi(t)$. For the case of the Riemann A(x), it was proved by Hardy that $\Xi(t)$ had an infinite number of real zeros using the standard discrete inversion symmetry of the underlying theta function, equivalent to the invariance of $A_1(x)$ defined above. Consideration of this symmetry of course implies that A(x) be defined not only on the interval $[1,\infty]$ but also on $[0,1]$. The basic character of this symmetry become clearer when the integral representation of $\Xi(t)$ is recast in the standard form of a cosine transformation after the change of variable $x = e^{2y}$:

$$\Xi(t) = 4\int_0^\infty dy A_I(e^{2y})\cos(ty) \tag{3.8}$$

The 1/x symmetry can now be understood as a reflection symmetry $y \to -y$ of $A_I(e^{2y})$, and if an A is analytic and a slowly varying quadratic around y=0, this is of course instrumental is producing an infinite number of real zeros in t. Generally if A(x) is chosen so that $A_I(x)$ is invariant under $x \to 1/x$, then the resulting $\Xi(t)$ is likely to have an infinite number of real zeros. This is so, for example, if $A_I(x)$ depends simply on the invariant combination x + 1/x, say like exp[-(x+1/x)] which leads to a $\Xi(t)$ involving Bessel functions of the third kind, as considered by Titchmarsh (and originally by Pólya, see [20]) when examining certain bona fide approximations to the Riemann $\Xi(t)$ function. It is quite useful to compare the properties of $\Xi(t)$ when choosing different A(x) in the integral representation class under study, so let us note for reference that

$$\Xi(t) = 2K_{it/2}(2) \text{ for } A_I(x) = Exp\left[-(x+\frac{1}{x})\right] \tag{3.9}$$

Moreover, different Bessel function examples of this type appear to have not only an infinite number of real only zeros but they also have an asymptotic density distribution similar to the one of the Riemann $\Xi(t)$. On the other hand, if the symmetry $x \to 1/x$ is not respected or analyticity at x=1 is violated, e.g. by choosing $A_I(x)$ to be simply exp(-x) then the resulting $\Xi(t)$ is a symmetrized incomplete gamma function, $\Xi(t) = \Gamma(it/2,1) + \Gamma(-it/2,1)$, which has a real zero at infinity only. This can be understood simply because $A_I(x)$ is not analytic at x=1 when extended by hand to x<1 to satisfy inversion symmetry. It is also worth noticing that if one only retains the first term in the theta function expansion of the Riemann A(x), then the resulting $\Xi(t)$ has only one real zero. As more terms are retained, more and more real zeros appear. The corresponding $\Xi(t)$ function simply corresponds to a sum of symmetrized incomplete gamma functions.

Another instructive example appears if $A_I(x)$ is chosen simply to be 1 with compact support in the x-interval [1, exp(2w)] for any w>0, with w=1 as a good representative example. Then we find

$$\Xi(t) = 4\sin(wt)/t \text{ for } A_I(x) = 1, x \in \left[1, e^{2w}\right] \tag{3.10}$$

which is of course a paradigm of the type of entire function that we would like to elucidate, i.e. functions with damped oscillations, an infinite number of real zeros, and no other zeros in the complex plane. If we renounce on analyticity at x=1 and choose $A_I(x)$ to be one at x=1 and decreasing linearly to zero on the same interval as above, i.e. not quite as spiked as the incomplete gamma function example discussed earlier, then for w=1 we find $\Xi(t)= 8(\sin(t/2)/t)^2$ and we have a limiting case when the real zeros (non-simple) are about to disappear completely. Another interesting intermediate case arises when $A_I(x)$ is chosen to be $\cos(\log(x)p/4)$ on the same interval as above for $w=1$, which leads to a quadratically decreasing $\Xi(t) = 2p \cos(t)/((p/2)^2 - t^2)$ having an infinite number of simple zeros.

In the present paper the main focus is to understand the properties of the Riemann $\Xi(t)$ function and its corresponding A(x) which is of exponential order 1 and type $p$. The analysis of other examples of A(x) are useful for understanding the larger picture of the considered integral representation. The examples of A(x) with compact support are of fundamental importance in understanding certain technical features of the polynomial approximation because of leading to more explicit expressions. The A(x) choices corresponding to $\Xi(t)$ being Bessel functions or incomplete gamma functions are more closely related to the Riemann case but are already quite complicated to analyze formally. It has been seen that it is sufficient to investigate A(x) functions of order 1, but it will turn out that it is useful to extend the analysis to cases with quite different exponential type, i.e. not only $p$ or 1 but much smaller ones as well.

A more comprehensive study of the A(x) function family also has to include some cleverly engineered examples of Dirichlet series, known in the literature for quite some time, which have properties comparable to the Riemann series but which nevertheless explicitly exhibit complex zeros of their corresponding $\Xi(t)$ function (see e.g. [24] and [8]). In this respect we would here like to note two different one-parameter families of $A_I(x)$ with this property. The first one, denoted $A^{(1)}$, is related to the Ramanujan tau function [7] and is explicitly given by

$$A^{(1)}(x,k) = \left[ x^{1/8} \exp(-p\sqrt{x}/12) \prod_{n=1}^{\infty} (1 - \exp(-2pn\sqrt{x})) \right]^k \quad (3.11)$$

where k is a positive integer different from 1, 2, 3, 4, 6, 8, 12 and 24. Such an A function is of exponential order 1/2, decreases monotonously like the Riemann A(x), and corresponds to a $\Xi(t)$ function which has an infinite number of real zeros, but which also has complex zeros outside the critical strip.

A second one- parameter $A_I(x)$ family of exponential order 2, denoted $A^{(2)}$, is given by [5]:

$$A^{(2)}(x,k) = \sum_{n=1}^{\infty} \left[ kA^{(2,1)}(x,n) - c(n)A^{(2,2)}(x,n) \right]$$
$$A^{(2,1)}(x,n) = (4p^2n^4x^4 - 6pn^2x^2)\exp(-pn^2x^2) \quad (3.12)$$
$$A^{(2,2)}(x,n) = (4p^2n^4x^4/25 - 6pn^2x^2/5)\exp(-pn^2x^2/5)$$

where $c(n)$ is the quadratic Dirichlet character modulo 5 ($c(1) = c(4) = 1, c(2) = c(3) = -1, c(5) = 0$) and $k \geq 4$. We note that this A function does not decrease monotonously like the Riemann A(x) does. The corresponding $\Xi(t)$ function has no real zeros and presumably an infinite number of complex zeros inside and outside the critical strip.

In summary, the above remarks serve to illustrate that the Riemann $\Xi(t)$ is part of a broad family of entire functions which may have a quite similar structure of real/complex zeros, as well as a very different one. When studying polynomial approximants to $\Xi(t)$ the key issue to be investigated below is how it is possible to arrange a most favorable setting for the approximants to have real roots only. If this can be done then of course different things may happen in the limit, some or all of the approximant roots may converge, or they may all disappear to infinity. But more importantly, there would be some hope that for specific choices of A(x), including the Riemann A(x), the $\Xi(t)$ could be shown to have no non-real zeros.

Starting from the convergent expansion of $\Xi(t)$ in symmetrized Pochhammer polynomials stated above, valid for any $b > 0$, we will now examine the polynomial approximants

$$\Xi_n(t, b) = \sum_{k=0}^{n} b_k(b) P_k^+(t/b) \qquad (3.13)$$

which have the property that $\Xi_{2n}(t, b)$ as well as $\Xi_{2n+1}(t, b)$ are even alternating polynomials of degree n in $t^2$, as it can easily be seen from the general properties of the polynomials $P_k^+(t)$. Let us note that we have the general decomposition

$$P_k(it) = P_k^+(t) + i P_k^-(t) \qquad (3.14)$$

where $P_k^+(t)$ is the even part of $P_k(it)$, as defined earlier, and $P_k^-(t)$ is the odd part of $P_k(it)$ divided by i. For reference let us write out explicitly the lowest k expressions:

$$P_0^+(t) = 1, P_1^+(t) = 1, P_2^+(t) = 1 - t^2/2, P_3^+(t) = 1 - t^2 \qquad (3.15)$$

$$P_0^-(t) = 0, P_1^-(t) = -t, P_2^-(t) = -3t/2, P_3^-(t) = -t(11 - t^2)/6 \qquad (3.16)$$

Using the basic recursion relation

$$P_{k+1}(it) = P_k(it)(1 - it/(k+1)) \qquad (3.17)$$

one finds the even/odd recursions

$$P_{k+1}^+(t) = P_k^+(t) + t P_k^-(t)/(k+1), P_{k+1}^-(t) = P_k^-(t) - t P_k^+(t)/(k+1) \qquad (3.18)$$

It is quite easy to see that $P_k^+(t)$ and $P_k^-(t)$ have real roots only. For example if $P_k^+(t) = 0$ then it follows when putting t = u + iv, with u and v real, that

$$|P_k(iu - v)| = |P_k(-iu + v)| \qquad (3.19)$$

and therefore

$$\prod_{j=1}^{k}((j+v)^2 + u^2) = \prod_{j=1}^{k}((j-v)^2 + u^2) \qquad (3.20)$$

But if v were non-zero, then the terms on the right hand side of the equation would all be bigger or all smaller than the corresponding terms on the other side. We thus conclude that v=0 and so the root t is real. Also we see that $P_k^+(t)$ and $P_k^-(t)$ cannot have any common root, because if they had one, then $P_k(it)$ would have a real root, which is not the case.

The above results also follow from the very structure of the recursions relations which imply generally that the polynomials $P_k^+(t)$ and $P_k^-(t)$ have interlacing distinct real roots and also form an interlacing sequence (see e.g. [13] for a comprehensive review of interlacing polynomials, as well as [20]). An alternative simple way of proving that $P_k^+(t)$ and $P_k^-(t)$ have interlacing distinct real roots is to apply the Hurwitz theorem for positive polynomials which have negative real roots only (these polynomials are also commonly called stable polynomials or Hurwitz polynomials). The theorem states that the even and odd part of such polynomials are interlacing polynomials of the same positive type. One can then use the one-to-one correspondence between positive polynomials with negative roots only and alternating polynomials with real roots only to prove the statement.

It is well-known that linear combinations with real coefficients (in fact both positive and negative ones) of interlacing polynomials with real roots only, produce polynomials with real roots only and which also interlace with the constituents. This circumstance is of course quite important for the purpose of investigating the real root properties of the polynomial approximants of $\Xi_n(t, \boldsymbol{b})$ introduced above. Indeed when examining specific cases numerically (including the Riemann approximants) one always finds that the lower order approximants $\Xi_n(t, \boldsymbol{b})$ have real roots only in t.

This situation should be contrasted with the case of a standard Taylor expansion of the Riemann $\Xi(t)$, as alluded to earlier, where partial sums of the cos(t log(x)/2) expansion are polynomials with complex roots only. Also other more standard attempts of expanding $\Xi(t)$ in a series of polynomials, such as Meixner-Pollaczek polynomials [16], show that complex polynomial roots are the rule rather than the exception. A popular remedy to such situations is to apply so-called multiplier sequences which transform a polynomial into another polynomial with better real root properties. Unfortunately these remedies are not generally applicable and often appear to be rather ad-hoc. The novel point in the approach of this paper is the presence of a continuous parameter β in the approximants which provides a more powerful degree of freedom for tailoring their real root properties.

Even if the real root structure is favorable when using approximants $\Xi_n(t, \boldsymbol{b})$ which are linear combinations of symmetrized Pochhammer polynomials $P_k^+(t)$, having distinct, real only and interlacing roots, there is still the complication that when n becomes sufficiently large then complex roots slowly make their appearance. This circumstance can be understood in a simple way. Let us denote the increasing sequence of squared roots of $P_k^+(t)$ by $r_{k,j}$, i.e. $r_{k,j} < r_{k,j+1}$ since the roots are distinct. The $P_k^+(t)$ sequence in k has interlacing roots, i.e. $r_{k+1,j} < r_{k,j} < r_{k+1,j+1}$. But unfortunately it turns out that the $P_k^+(t)$ sequence is not totally interlacing in the sense of having all the roots with number j belonging to disjoint intervals. For example, the first root squared $r_{2,1} = 2$ is smaller than the second roots squared $r_{k,2}$ only until k=23.

If the $P_k^+(t)$'s had been a totally interlacing sequence then there would have been quite good reason to expect the approximants $\Xi_n(t, b)$ to be polynomials with real roots only. But the interlacing property is only valid in a sequential neighborhood, albeit in a rather broad one, and therefore the approximants $\Xi_n(t, b)$ would logically be expected to feature some complex roots from a certain n and beyond. Fortunately we will see in the next section that, even without the totally interlacing property, it is still possible to infer something about the real root characteristics of $\Xi_n(t, b)$ due to a particular feature of the coefficients $b_k(b)$.

Before addressing this feature, let us for completeness just note that if a certain (presumably quite large) $b$ would exist such that the polynomials $\Xi_n(t, b)$ for sufficiently large n had only real roots, then we could already at this point jump to our final conclusion about the zeros of the function $\Xi(t)$ to which $\Xi_n(t, b)$ is then certain to converge. In the following sections we will therefore assume that such a fixed beta does not exist, so that the central issue will rather be about the properties of n-dependent beta-sequences $b_n$. Section 4 will examine the characteristics of $b_n$ with respect to the onset of a certain real root regime for the polynomial approximants, while section 5 will examine the conditions on $b_n$ assuring convergence of the approximants to the $\Xi(t)$ function.

## 4. Existence of beta-sequences preserving real roots

In the expression for the approximant $\Xi_n(t, b)$ defined above, we have a particular positive linear combination of $P_k^+(t)$'s and one may wonder if there happens to be any generic, finite summation formula for the $P_k^+(t)$'s which resembles this expression. There is an affirmative answer to this question which can be found by defining the sum $S_n(it)$ by

$$S_n(it) = \sum_{k=0}^{n} \frac{1}{k+1} P_k(it) \qquad (4.1)$$

Using the basic recursion formula rewritten as $itP_k(it)/(k+1) = P_k(it) - P_{k+1}(it)$ and summing over k we obtain

$$itS_n(it) = 1 - P_{n+1}(it) \qquad (4.2)$$

If we now decompose in even and odd parts $S_n(it) = S_n^+(t) + iS_n^-(t)$, then one easily finds

$$tS_n^+(t) = -P_{n+1}^-(t), -tS_n^-(t) = 1 - P_{n+1}^+(t) \qquad (4.3)$$

In conclusion, we get the explicit summation formula

$$S_n^+(t) = \sum_{k=0}^{n} \frac{1}{k+1} P_k^+(t) = -P_{n+1}^-(t)/t \qquad (4.4)$$

which states that a harmonic sum of symmetrized Pochhammer polynomials is simply proportional to an anti-symmetric Pochhammer polynomial $P_k^-(t)$. The important result follows directly that the sum $S_n^+(t)$ is a polynomial with distinct real only roots.

Let us now see how this result impacts on our understanding of the real root characteristics of the approximant $\Xi_n(t,b)$. With a change of variables $y=x^{-b/2}$ we rewrite the coefficient $b_k(b)$ as

$$b_k(b) = \frac{4}{b} \int_0^1 dy\, A_I(y^{-2/b})(1-y)^k \tag{4.5}$$

We notice that A now appears directly in terms of the expression $A_I$, just like in the cosine representation of $\Xi(t)$. The limit $x \to \infty$ corresponds to the limit $y \to 0$ which will be critical for understanding the asymptotic properties of the series.

The surprisingly simple key observation is now that when $b$ becomes large then the integrand term depending on $b$ can be replaced by a constant and we simply get

$$b_k(b) \approx \frac{4A_I(1)}{b} \int_0^1 dy\,(1-y)^k = \frac{4A_I(1)}{b(k+1)} \tag{4.6}$$

For large $b$ we therefore find that

$$\Xi_n(t,b) \approx \frac{4A_I(1)}{b} \sum_{k=0}^n \frac{1}{k+1} P_k^+(t/b) = \frac{-4A_I(1)}{t} P_{n+1}^-(t/b) \tag{4.7}$$

which is an even polynomial with distinct real only roots in t. This observation of $\Xi_n(t,b)$ having distinct real roots only in the asymptotic scaling limit $b \to \infty$ will be seen below to have important consequences.

It is quite easy to verify numerically all the above observations concerning real polynomial roots using a standard PC software package, like Wolfram Mathematica, for the case of the Riemann A(x), as well as for the other admissible A(x) discussed above. The important novel feature to be examined in detail below is the characterization of the precise range of $b$ values for which the polynomial approximants have real roots only. Taking the most interesting example of the Riemann A(x), one indeed finds that the approximants $\Xi_n(t,b)$ for low n (not less than 4), starting with a rather small $b$ around 1-3, are polynomials with distinct real only roots (in a t range starting from around 8). As one increases n, at some point a pair of complex roots generally appear in the lower root range, and one may then try to increase $b$ slightly until the complex root pair coalesces into a real double root, and all roots again become real, and remain real for all larger $b$ values.

In a separate paper [10], a detailed numerical analysis, including graphs, of several different A(x) cases will be presented based on the use of PC as well as mainframe computing resources. This kind of analysis can be carried on with relative ease to n of the order of millions, but here we will mainly limit the discussion to how one can extract the salient root features of the polynomial approximants for relatively low n using generally available computing tools. Following the above iterative numerical procedure in the range of n from 10 to 200, we start by evaluating $b_{10} \sim 3.23$ and $t_{10} \sim 10.63$ and proceed by incrementing n by one. In this manner, we find that the onset of the real root regime, i.e. the maximal connected range of beta values for which there are no complex roots for a chosen n, for the Riemann A(x) emerges for beta values according to the following approximate fit (with an $R^2$ of 0.9999998):

$$b_n \approx 5.58(\log(n+1))^{0.66} - 6.76 \quad (4.8)$$

If we analyze the Bessel function example $K_{it/2}(2)$, then we may start by evaluating $b_{10} \sim 2.12$ and $t_{10} \sim 7.07$ and, proceeding as before, we find the following fit to the onset of the real root regime

$$b_n \approx 3.95(\log(n+1))^{0.63} - 4.74 \quad (4.9)$$

For these A's which correspond to $\Xi(t)$ functions having infinitely many real roots, we note that the $b_n$ sequence is monotonously increasing for all the (low) n values considered. Similarly we find that the positive double roots $t_n$, the coalescing point for complex root pairs defining the onset of the real root regime, also form an increasing sequence located in the range 10-13, with a quite slow approach towards the lowest zero of the corresponding $\Xi(t)$ (i.e. $\sim 14.1$ for the Riemann $\Xi(t)$). For these representative $\Xi(t)$ cases, including the Riemann one, we thus find that one seems to have sub-logarithmic growth rate $(\log(n+1))^c$ with $c \sim 0.6$, with a rather smooth, increasing behavior of $b_n$ and a very close numerical logarithmic power fit of $b_n$ in the range of n between 10 and 200.

However, there are other examples where the numerical study in the same n-range of 10-100 shows a clearly differentiated behavior of the $b_n$ sequence defining the onset of the real root regime. For the case of the A(x) of a simple exponential form, corresponding to a $\Xi(t)$ which is a symmetrized incomplete gamma function with no real roots, we find the growth rate $c \sim 0.86$. Here the polynomial roots $t_n$ again form a monotonous, slowly increasing sequence.

For the sin(t)/t example corresponding to a step function $A_1(x)$, we start the evaluation with the lowest possible n = 4 and find $b_4 \sim 0.07$ and $t_4 \sim 3.46$. Proceeding as above, we find the following fit to the onset of the real root regime

$$b_n \approx 1.02(\log(n+1))^{0.99} - 1.60 \quad (4.10)$$

In this case we also note some narrow patches around some particular n-values which are characterized by a slight non-monotonous behavior of $b_n$. This phenomenon is related to the peculiar circumstance that at these n-values, the position of the real double root in $t_n$ suddenly jumps to a higher value located in the same interval of the zeros of sin(t), i.e. the integers times π. For example, at n=6 one observes that the double-root of the approximant jumps from 3.47 to 5.17. After that, $t_n$ increases more slowly for a while, and then at n=17 it jumps again from 6.69 to 7.93. This phenomenon is repeated ad infinitum as will be discussed later.

A little before the n-value corresponding to the t-jump, then $b_n$ is observed to decrease very slightly. Subsequently it increases monotonously until reaching the next jump location. The fact that the $b_n$ behavior is slightly irregular is responsible for the numerical fit for the sin(t)/t case not being a good as for the Bessel and Riemann case. But as a matter of fact, it will be argued below that for the sin(t)/t case one would expect to find an exact $b_n$ asymptotic growth rate of $(\log(n+1))^c$ with c=1, so that here the low n behaviour is nevertheless a reliable guide to what happens at very high n. On the other hand, although the situation for the Bessel and Riemann cases looks simpler because of the monotonous smooth $b_n$ behaviour, these cases are rather deceptive. In fact, the $b_n$ sequence continues to behave smoothly for quite a while [10], with no t-jump effect to seen below n equal to two million! However, we will show below that eventually there must be t-jumps like in the sin(t)/t case.

The $A^{(1)}(x,k)$ case for k=5, for example, related to the Ramanujan tau function has some points in common with the sin(t)/t case. In the range of n<200, we may start by calculating $b_{10} \sim 1.006$ and $t_{10} \sim 3.47$ and proceed to find that the $b_n$ sequence is increasing smoothly, except for n = 66 where it jumps downwards, and where the $t_n$ root is seen to jump upwards within an interval of the lowest $\Xi(t)$ zeros located at 3.83 and 6.41. The best fit to the growth rate of the $\beta_n$ sequence shows a supra-logarithmic behaviour of $(\log(n+1))^{1.30}$. Finally for the $A^{(2)}(x,k)$ case for k=5, for example, with no real zeros as mentioned earlier, we find no t-jumps and a rather smooth monotonous behavior of the $b_n$ sequence with a supra-logarithmic growth rate of $(\log(n+1))^{1.32}$.

It should be emphasized that all of the above numerical statement are based on an iterative procedure which can be carried through to an arbitrary numerical precision. The rate of growth of the $b_n$ sequence will be subject to further discussion below, but clearly the presented growth rate values c are only indicative because of the small range of n-values studied which leaves some room for alternative fits. It is noteworthy, however, that the study of $b_n$ sequences for different types of A(x) already for quite low n-values allows to distinguish rather clearly four different types of behaviour. There is a class of behaviour with c<1 of distinct sub-logarithmic beta-sequence growth, including a sub-class of $\Xi(t)$ cases with infinitely many zeros, like the Riemann case, with a growth rate c ~ 0.6 for n<200, and a sub-class with cases with A(x) corresponding to $\Xi(t)$ with no real zeros, or just a finite number of real zeros, for which c ~ 0.85 for n<200. Finally there is class of behaviour with c≥1 of logarithmic or supra-logarithmic growth rate, including a sub-class of functions like the sin(t)/t case and other cases of compact support for A(x) with c=1, and sub-class of some more exceptional $\Xi(t)$ cases corresponding to special Dirichlet series with c>1.

This initial rough differentiation of A(x) function classes may be pursued in two directions. Firstly, one may subject a number of representative A(x) to more detailed numerical scrutiny so a to better ascertain the $b_n$ growth rate for very high n [10], and in particular to discuss how confident one may be about the kind of approximate analytical fits discussed above. The log-power type of fit is clearly only one of several possibilities which may be considered for the finite range analysis of a numerical approach. In fact, there are other candidates, like the type of sublog(n) and sublogxl(n) functions to be discussed later, which may provide equally good fits for low n, and which moreover are much more likely to be close to the ultimate asymptotic behaviour of $b_n$. Secondly, one may develop a more formal analytical approach in terms of specific difference equations for the $b_n$ and $t_n$ sequences. This will be done below in section 6, 7 and 8.

Before going into this, we may summarize the situation of our present analytical understanding of the β dependence of the approximants by first noting it has quite generally been established that in the asymptotic scaling limit $b \to \infty$ then the polynomial approximant $\Xi_n(t,b)$ has distinct real only roots. If we examine for any given n what happens in root space when we decrease the continuous parameter $b$ from the asymptotic range down to the smaller $b$ range (see [21] for a general discussion of root space analysis), then we note that the roots of $\Xi_n(t,b)$ must follow continuous trajectories in the complex plane, and clearly distinct real roots will remain distinct real roots for quite some time. What typically happens from some point on is that two low lying real roots coalesce and become a real double root. Subsequently this root may separate into to a pair of complex roots and we then enter a mixed real/complex root regime.

The above remarks based on standard root space analysis lead us to conclude directly that there must exist increasing beta-sequences $b_n \to \infty$ such that the polynomial approximant $\Xi_n(t,b_n)$ has distinct real only roots for every n. Moreover, for each given A(x), we can state that there exists a unique minimal beta-sequence $b_{min,n}$, such that $\Xi_n(t,b_{min,n})$ has real only roots, which is simply defined by choosing the first $b$ for which a real double (or multiple) root of $\Xi_n(t,b)$ appears, often, but not necessarily, on the boundary to the mixed root regime, when moving down from the asymptotic limit. Since $b_{min,n}$ is defined as the infimum of a connected semi-infinite range, any beta-sequence above this minimum sequence, i.e. $b_n > b_{min,n}$, will produce a $\Xi_n(t,b_n)$ having distinct real only roots for every n.

The empirical beta sequence data and approximate fits studied above are exactly the low n manifestation of this unique minimal beta-sequence characteristic of each individual A(x). Before analyzing the more detailed characteristics of minimal beta-sequences, we will however first discuss under which conditions the polynomial approximants involving an increasing $b_n$ sequence actually converge to $\Xi(t)$.

## 5. Convergence of the approximants

The uniform convergence of $\Xi_n(t,b)$ to $\Xi(t)$ on any compact subset of the complex t-plane has already been established for fixed $b$, but we are now interested in analyzing whether the approximant $\Xi_n(t,b_n)$ might still converge to $\Xi(t)$ for a certain range of not too quickly growing beta-sequences $b_n \to \infty$ as $n \to \infty$. If this is so, then of course it would also be interesting to check if this range of beta-sequences happens to include sequences with $b_n \geq b_{min,n}$ for sufficiently large n, so that $\Xi_n(t,b_n)$ would have real roots only, as well as having the property of converging to $\Xi(t)$.

There is in fact good reason to believe that $\Xi_n(t,b_n)$ might actually still converge when moving from the fixed $b$ case through a range of slowly increasing beta sequences to a certain maximal limit case.

Let us first, more intuitively, consider the situation for t = 0 when $P_k^+(0) = 1$ so that, as previously found, we have for large $b_n$

$$\Xi_n(0, b_n) \approx \frac{4A_I(1)}{b_n} \sum_{k=0}^{n} \frac{1}{k+1} \approx \frac{4A_I(1)}{b_n}(\log(n+1) + C) \tag{5.1}$$

where C is Euler's constant. If there is convergence to the bona fide limit value $\Xi(0)$ then we would here intuitively expect that $b_n \sim B \log(n)$ and that $B = 4A_I(1)/\Xi(0)$. If one had chosen B to be larger than this constant, then the limit value for $n \to \infty$ would clearly be too small, and one could then argue that B were too big to neglect the remainder term. On the other hand, if B were chosen smaller that this constant, then the limit value would be too large, but then one could argue that we were anyway getting closer to the fixed $b$ case where a more careful analysis of the sum would show convergence to the correct value.

The simple intuition turns out to be approximately right, but clearly a more careful convergence analysis of the $n \to \infty$ limit is needed. We can prove that $\Xi_n(t, b_n)$ converges to $\Xi(t)$ if it is possible to show that the remainder term $R_n(t, b_n) = \Xi(t) - \Xi_n(t, b_n)$, or explicitly

$$R_n(t, b_n) \equiv \Xi(t) - \Xi_n(t, b_n) = \sum_{k=n+1}^{\infty} b_k(b_n) P_k^+(t/b_n) \tag{5.2}$$

goes to zero as $n \to \infty$ on any compact subset of the complex plane.

Looking first at the case t = 0, as above when $P_k^+(0) = 1$, we carry out the summation over k in the exact expression for $R_n(t, b_n)$ to get

$$R_n(0, b_n) = \frac{4}{b_n} \int_0^1 dy\, A_I(y^{-2/b_n})(1-y)^{n+1}/y \tag{5.3}$$

The important part of the integration is for y small where the $A_I$ factor in the integrand provides an exponential damping which is getting weaker for large n. The structure of the large n limit become a little clearer if we consider the specific example of $\Xi(t) = 4\sin(wt)/t$ for which $A_I(x) = 1$ on the sub-interval $[\exp(-wb_n), 1]$:

$$R_n(0, b_n) = \frac{4}{b_n} \int_{\exp(-wb_n)}^{1} dy(1-y)^{n+1}/y \tag{5.4}$$

The large n asymptotics of this integral can be found, for example, by a straightforward application of the Laplace method after a partial integration. For the integral to go to zero for large n, the lower limit of the y-integral cannot be too small, which means that $\beta_n$ cannot grow too fast with n. If we assume that it is of the form $b_n = B \log(n)$ and set a =wB then the lower limit of the integral I(n,a) is $n^{-a}$ and we find for large n that I(n,a) for $a \leq 1$, while $I(n,a) \sim \log(n)$ for a>1. Therefore we see that the full remainder term for $a \leq 1$ behaves like

$$R_n(0, b_n) \approx \frac{n^{a-1}}{\log(n)} \tag{5.5}$$

for large n, while for a>1 it does not go to zero.

In particular for the critical value a=1 the convergence of $\Xi_n(0, b_n)$ to $\Xi(0)$ is logarithmic. The qualitative behavior of the integral is clearly determined by the power factor $(1-y)^{n+1} \sim e^{-1}$ for $y^{-1} = n$ and n large, which will also be seen below to be an interesting point in the convergence analysis of the general $A_I$ case. Finally let us note that if $b_n$ grows slower than $1/w \log(n)$ then there is of course also convergence, all the way down to the fixed $b$ situation.

For other choices of $A_I$ with non-compact support and genuine exponential decrease, the large n analysis of the $R_n(0, b_n)$ integral is slightly more complex, but in essence the exponential damping acts as a small y cut-off and a similar convergence limit feature as in the $A_I = 1$ case emerges. The convergence analysis of the remainder term can be carried through by establishing asymptotic properties of upper bounds, so it is sufficient to examine the generic case $A_I(x) = \exp(-ax)$ if $A_I(x)$ is of exponential order one, e.g. for the case of the Riemann $A(x)$. In fact, the analysis is quite identical for a general normal exponential order of $A_I(x) = \exp(-ax^b)$ with b>0, as one might also expect from the already established power transformation equivalence relation in the $A(x)$ family.

For the exponential case, the significant difference with the compact support case is that the critical growth rate of $b_n$ for convergence cannot be logarithmic. In fact, convergence can in this case only be assured for a sub-logarithmic growth rate

$$b_n = o(\log(n)) \tag{5.6}$$

A numerical study of the convergence properties of the remainder term may superficially indicate that logarithmic growth of $b_n$ is sufficient to assure convergence, but this turns out be deceptive because of very slow log(log(n)) features. To prove the above statement on the convergence rate, we will proceed in two steps.

As a first step, one may note that the analytical approach to the y-integral of the remainder term given above naturally suggests a specific growth rate of $b_n$ which is sufficient for assuring convergence, even it may not be an optimal growth rate. When applying the Laplace method to the integrand of the remainder term with $A_I(x) = \exp(-ax)$, which has a nicely skewed bell-shaped y-dependence, we have to analyze an exponential term of the form

$$f(y) = -\log(y) + (n+1)\log(1-y) - ay^{-e(n)} \tag{5.7}$$

where for short we have set $e(n) = 2/b_n$. Taking the derivative of this expression, one is led to look for approximate solutions to a transcendental equation containing an expression

$$y^{-1} = (1/e(n))^{1/e(n)} \qquad (5.8)$$

From the convergence discussion of the $A_I = 1$ case we have seen that an extremal value of $y^{-1}$ should preferably be of order n. But since $\log(n)^{\log(n)} = n^{\log(\log(n))}$, it is in this case warranted to try to select a more appropriate $b_n$ expression which has a slightly slower growth rate than log(n). For lack of a generally accepted terminology, we will here denote by sublog(n) the function which is a (principal) solution of the transcendental equation

$$(\text{sublog}(n))^{\text{sublog}(n)} = n \qquad (5.9)$$

This function is an example of a so-called tetration (or hyper-4, iterated exponentiation) which has also sometimes been referred to as a super-square-root or a super-log. It can be expressed in terms of the more generally known Lambert W-function through the relation

$$\text{sublog}(n) = \frac{\log(n)}{W(\log(n))} \qquad (5.10)$$

which can easily be verified by using the defining equation W(n) exp(W(n)) = n. A good reason for choosing to denote the function by sublog(n) is that it grows only very slightly slower than log(n), as expressed by the following growth properties

$$\text{sublog}(n) = o(\log(n)), [\log(n)]^d = o(\text{sublog}(n)) \text{ for } d < 1 \qquad (5.11)$$

For completeness, let us note the useful definition including an additional parameter $r > 0$

$$\text{sublog}(n, r) \equiv r\,\text{sublog}(n^{1/r}) \qquad (5.12)$$

which fulfils the implicit equation encountered in many situations

$$\log(n) = \text{sublog}(n, r) \log(\text{sublog}(n, r)/r) \qquad (5.13)$$

In summary, if we suppose that the beta-sequence is precisely of the form

$$b_n = 2\,\text{sublog}(n) \qquad (5.14)$$

then the Laplace method can be applied to show that the remainder term for an $A_I$ of exponential order one and type a>0 has the large n form

$$R_n(0, b_n) \approx \frac{2e^{-a-a\,\text{sublog}(n)}}{\text{sublog}(n)} \qquad (5.15)$$

which clearly goes to zero for $n \to \infty$ with a good margin. Of course, if the sequence $b_n$ has a slower growth rate than this, then the remainder term will also go to zero for $n \to \infty$. It is also quite easy to see that in the more general case of $A_I$ being of exponential order b, i.e. $A_I(x) = \exp(-ax^b)$, then we can simply choose

$$b_n = 2b \operatorname{sublog}(n) \tag{5.16}$$

to assure that the remainder term goes to zero for large n. As already remarked, the choice of $\varepsilon(n) = 2/b_n = 1/\operatorname{sublog}(n)$ is sufficient to assure convergence, but there is good reason to suspect that there may also be convergence beyond this rate. The interesting feature of the sublog(n) function is that it represents as simple dividing line between the slower sub-logarithmic growth rate of the type $(\log(n+1))^c$ with c<1 and the faster sub-logarithmic growth rates of the type sublogxl(n) to be considered later. As a second step in the argument, let us therefore consider the integral

$$I_n = \int_0^1 dy \exp(-ay^{-e(n)})(1-y)^{n+1}/y \tag{5.17}$$

and assume that $\operatorname{sublog}(n) \ll 1/e(n) \ll \log(n)$. Splitting the integral in two part $I_n = I_n^1 + I_n^2$ where the first term corresponds to the integration interval [0,1/n] and the second term to the integration interval [1/n,1], it is straightforward to apply the Laplace method to find that $I_n^2 \sim \exp(-an^{e(n)})$ which goes to zero for $n \to \infty$ since $1/\varepsilon(n) \ll \log(n)$ and therefore $n^{e(n)} \to \infty$. The first term $I_n^1$ can be majorized by

$$I_n^1 \leq \int_0^{1/n} dy \exp(-ay^{-e(n)})/y \tag{5.18}$$

and after a change of variables $y = u^{-1/e(n)}$ it is easily seen that $e(n) I_n^1 \to 0$ for $n \to \infty$ because $n^{e(n)} \to \infty$. In conclusion, the remainder term $R_n(0, b_n)$ will go to zero if we have $b_n = o(\log(n))$. Also one can easily convince oneself that there is not convergence for $b_n \sim \log(n)$.

It is not too difficult to extend the above results for $\Xi_n(t, b_n)$ at t = 0 and the corresponding remainder term to any compact subset of the complex plane. Let us first recall the previously used uniform estimate for large k:

$$|P_k(it/b)| < C_1 k^{\operatorname{Im}(t)/b} \tag{5.19}$$

valid on any compact subset. If $|\operatorname{Im}(t)| \leq M$ (e.g. M=1/2 for the critical strip), we therefore have the bound

$$|P_k^+(t/b)| < C_1 k^{M/b} \tag{5.20}$$

Consequently we find that

$$|R_n(t,b_n)| < \frac{4C_1}{b_n} \int_0^1 dy\, A_I(y^{-2/b_n})(1-y)^{n+1} \Phi(1-y, -M/b_n, n+1) \tag{5.21}$$

where the Hurwitz-Lerch function $\Phi$ is given by

$$\Phi(1-y, -M/b_n, n+1) = \sum_{k=0}^{\infty} (k+n+1)^{M/b_n} (1-y)^k \tag{5.22}$$

Since the parameter $M/b_n \to 0$ for large n we note that $\Phi(1-y, -M/b_n, 2n+1) \sim 1/y$ and a similar convergence analysis as for t = 0 can be applied to show that the remainder term goes to zero on any compact subset provided that $b_n$ does not grow faster than discussed as above. Let us briefly examine the possible adverse convergence effect of having M>0 in the two cases of $A_I(x)=1$ and $A_I(x) = \exp(-ax)$.

In the case of $A_I(x)=1$ where the maximal growth rate is $b_n = B \log(n)$, then we can estimate the significance of the M dependent term by examining the following approximation

$$(k+2n+1)^{M/b_n} = \exp\left[\frac{M}{b_n}\log(k+n+1)\right] \approx \exp\left[\frac{M}{B}\right]\left[\exp\left[\frac{M}{Bn\log(n)}\right]\right]^k \tag{5.23}$$

The effect of M on the remainder term can here be seen in an overall factor exp(M/B) as well as in a modification of the $\Phi$ power series which is effective for values of y of the order of M/(Bn log(n)). But for large n this is smaller that the lower integral limit of $n^{-a}$ when $a \leq 1$, so we conclude that the convergence properties are unchanged for M>0. This can also be verified by numerical evaluation of the integral involving the exact $\Phi$ factor.

The discussion for the case of $A_I(x) = \exp(-ax)$, where the beta sequence is of the form $b_n$ =2sublog(n), is quite similar. The M dependent factor can be estimated in the same way, but in this case it will give rise to an overall factor which grows like exp(log(sublog(n))). This is however much slower than the previously found decreasing factor of the remainder term. Also the modification of the $\Phi$ power series is effective only for values of y of the order of M/(2n sublog(n)) and therefore does not modify the convergence analysis. Again numerical evaluation of the integral with the exact $\Phi$ confirms this analysis.

Thus we conclude that in all $A_I(x)$ cases, including the Riemann A(x) case, the remainder term goes to zero on any compact subset of the complex plane. For the $A_I(x) =1$ case there is a maximal beta sequence of the form $b_{max,n} = 1/w\, \log(n)$, while in the case $A_I(x) = \exp(-ax)$ of exponential order 1, we must simply require $b_n$ =o(log(n)). In the case of the $A_I(x) =1$, we of course still have the freedom to add a constant to the maximal rate so that we may state that convergence is assured for any $b_n$ choice of the form $B_0 + b_{max,n}$ where $B_0$ is an arbitrarily (large) constant. In summary, we have established that for any A(x) belonging to the class of admissible functions, then $\Xi_n(t,b_n)$ converges uniformly to $\Xi(t)$ if the growth rate of $b_n$ is either smaller than or equal to a certain logarithmic beta sequence $b_{\max,n}$, or more restrictively $b_n = o(\log(n))$ for the exponential cases.

## 6. Asymptotic properties of minimal beta-sequences

The important question is now whether, for any given A(x), it is possible to develop a more precise analytical understanding of the n-dependence of the unique minimal beta sequence $b_n = b_{\min,n}$ associated with A(x), as well as of the associated value of the polynomial double root $t_n$. While the established fact that the polynomial approximant roots are all real and distinct for large $b$ is sufficient to assure the existence of such a minimal beta-sequence, it appears to be quite unfeasible to find any exact solution for this. Fortunately, it is only warranted to understand the asymptotic behaviour of the minimal $b_n$ sequence, but even that is far from straightforward.

We will show below how it is possible to find a simple recursive scheme for the minimal beta and t-sequences which is instrumental for both performing a detailed numerical analysis by computer as well as for formally disentangling the asymptotic behaviour of the sequences for large n. The basic idea is to mimic the iterative procedure described earlier which provided the exact numerical values (i.e. to within the chosen numerical precision) for $b_n$ and $t_n$ in the low n range of 4-100. Let us therefore suppose that for a given n, we have already found the values $b_n$ and $t_n$ corresponding to the situation where the polynomial $\Xi_n(t,b)$ has a real double root, such that for all $b > b_n$ the roots are all real and distinct while for $b < b_n$ complex roots start to appear. Consequently we have that both $\Xi_n(t,b)$ and its first derivative with respect to t, denoted by $\Xi'_n(t,b)$, are both zero at this point:

$$\Xi_n(t_n, b_n) = 0 \text{ and } \Xi'_n(t_n, b_n) = 0 \tag{6.1}$$

If we now increment n by one, we have

$$\Xi_{n+1}(t,b) = \Xi_n(t,b) + b_{n+1}(b) P^+_{n+1}(t/b) \tag{6.2}$$

and the values of $b_n$ and $t_n$ will also have to be incremented

$$b_{n+1} = b_n + \Delta b_n, \quad t_{n+1} = t_n + \Delta t_n \tag{6.3}$$

so as to satisfy the new double root condition

$$\Xi_{n+1}(t_{n+1}, b_{n+1}) = 0 \text{ and } \Xi'_{n+1}(t_{n+1}, b_{n+1}) = 0 \tag{6.4}$$

To first order in the increments one thus finds the following difference equations

$$\Delta b_n = -\frac{\Xi_{n+1}(t_n, b_n)}{\frac{\partial}{\partial b}\Xi_{n+1}(t_n, b_n)} \tag{6.5}$$

$$\Delta t_n = -\frac{\Xi'_{n+1}(t_n, b_{n+1})}{\Xi''_{n+1}(t_n, b_n)} = -\frac{\Xi_{n+1}(t_n, b_n)}{\Xi''_{n+1}(t_n, b_n)} \left[ \frac{\Xi'_{n+1}(t_n, b_n)}{\Xi_{n+1}(t_n, b_n)} - \frac{\frac{\partial}{\partial b}\Xi'_{n+1}(t_n, b_n)}{\frac{\partial}{\partial b}\Xi_{n+1}(t_n, b_n)} \right] \quad (6.6)$$

These equations simply encompass the classical Newton root finding method, and if they are iterated in the indicated order, keeping the polynomial order fixed, they provide a straightforward numerical approach for finding the double root values to any desired precision. For the numerical data mentioned earlier in the n-range between 10 and 100, the first order approximation is typically correct to within 3% and the Newton type iteration therefore converges rather quickly. The application of this approach and its practical adaptation to a mainframe computing platform is presented in a separate paper [10].

When n becomes very large then the polynomial increments are very small and therefore the first order approximation becomes better and better. Consequently, the formal asymptotic analysis of the beta and t-sequence can simply be applied using the above equations without iteration, in a concise form explicitly implementing the double root conditions $\Xi_n(t_n, b_n) = 0$ and $\Xi'_n(t_n, b_n) = 0$:

$$\Delta b_n = -b_{n+1}(b_n) Q_{n+1}(t_n, b_n) \quad (6.7)$$

$$\Delta t_n = -\frac{b_{n+1}(b_n) P^+_{n+1}(t_n/b_n)}{\Xi''_{n+1}(t_n, b_n)} \frac{\partial}{\partial t} \log(Q_{n+1}(t_n, b_n)) \quad (6.8)$$

$$Q_{n+1}(t_n, b_n) = \frac{P^+_{n+1}(t_n/b_n)}{\frac{\partial}{\partial b}\Xi_{n+1}(t_n, b_n)} \quad (6.9)$$

Below these equations will simply be referred to as the beta and t-increment equations. Here we have introduced explicitly the ratio $Q_{n+1}(t_n, b_n)$, the numerator and denominator of which will be seen to play an important role in the large n asymptotic analysis. Let us now try analyze in more detail the asymptotic properties of the individual terms of these difference equations. Firstly, the symmetrized Pochhammer polynomial $P^+_{n+1}(u_n)$, where have put $u_n = t_n/b_n$, can be approximated for large n by using the asymptotics of gamma function to write

$$P_n(it) \approx n^{-it} / \Gamma(1 - it) \quad (6.10)$$

Consequently

$$P^+_n(u) \approx \frac{\cos(u \log(n) + j(u))}{r(u)} \quad (6.11)$$

where we have set

$$r(u) = Abs(\Gamma(1 - iu)) = \left(\frac{pu}{\sinh(pu)}\right)^{1/2}, \; j(u) = Arg(\Gamma(1 - iu)) \approx Cu \quad (6.12)$$

and the last relation is valid for small u, with C being the Euler constant. Thus $P^{+}_{n+1}(u_n)$ behaves asymptotically in an oscillatory manner which depends on the limit value of $u_n$. Secondly, we recall that the coefficient $b_{n+1}(b)$ can be found by calculating an integral $B_n(b)$:

$$b_n(b) = \frac{B_n(b)}{b}, B_n(b) = 4\int_0^1 dy A_I(y^{-2/b})(1-y)^n \qquad (6.13)$$

where for convenience we have extracted the $1/b$ factor. The large n asymptotics of $B_n(b)$ can be evaluated by the same methods used for estimating the remainder term in the convergence analysis above. Thirdly, we note the appearance in the difference equations of the beta-derivative of $\Xi_{n+1}(t_n, b_n)$. This term is quite simple to include in a numerical analysis, but in a formal asymptotic analysis it appears superficially to be difficult to relate it to something simpler.

However, surprisingly it turns out that this derivative term can be related to another simple derivative function. To see this, let us first rewrite the infinite Pochhammer expansion of $\Xi(t)$ as

$$b\Xi(t) = \sum_{k=0}^{\infty} B_k(b) P_k^{+}(t/b) \qquad (6.14)$$

If one differentiates this expression with respect to $b$, as well as with respect to t, then one finds

$$\Xi^{*}(t) \equiv \Xi(t) + t\Xi'(t) = \sum_{k=0}^{\infty} B_k'(b) P_k^{+}(t/b) \qquad (6.15)$$

where $\Xi'(t)$ and $B_n'(b)$ denote derivatives with respect to the relevant variable, which hints that a $b$-scale transformation is somehow related to a t-scale transformation. That this is so should come as no surprise, since a $b$-scale transformation implies an x-power transformation of A(x), which is equivalent to a t-scale transformation of Ξ(t), as noted earlier. We can make this relationship more precise by evaluating the beta-derivative term appearing in the increment equations through the formula

$$b_n \tfrac{\partial}{\partial b} \Xi_{n+1}(t_n, b_n) = -[\Xi_{n+1}(t_n, b_n) + t\Xi'_{n+1}(t_n, b_n)] + \sum_{k=0}^{n+1} B'_k(b_n) P_k^{+}(t_n/b_n) \qquad (6.16)$$

Using the double root condition $\Xi_n(t_n, b_n) = 0$ and $\Xi'_n(t_n, b_n) = 0$, the first term of the right hand side evaluates to a combination of Pochhammer polynomials times a b-coefficient:

$$-b_{n+1}(b_n)\left(P^{+}_{n+1}(u_n) + u_n P^{+'}_{n+1}(u_n)\right) \qquad (6.17)$$

which asymptotically will be seen to be very small. On the other hand, the second term of the right hand side will for large n be quite close to $\Xi^{*}(t_n)$, as defined above, supposing of course that the beta-growth rate is slow enough to assure convergence of the approximant. We therefore conclude that the large n asymptotics of the beta-derivative term is closely related to the large t asymptotics of $\Xi^*(t)$ since $t_n \to \infty$ as $n \to \infty$. This particular, somewhat surprising, asymptotic feature

$$b_n \tfrac{\partial}{\partial b} \Xi_{n+1}(t_n, b_n) \approx \Xi^*(t_n) \tag{6.18}$$

is indeed confirmed, with a slight modification, by the numerical analysis as will be seen below.

In summary, we now have a reasonably good grasp on how to conduct the asymptotic analysis of all the different terms entering into the difference equations for the $\beta_n$ and $t_n$ sequences. Let us therefore start with analyzing the properties of the minimal beta sequence for the simplest reference case of $\Xi(t)$ =4sin(t)/t, when choosing $w=1$ in the more general case. For this case we have the explicit formula

$$b_n(b) = \frac{4}{b(n+1)}(1-e^{-b})^{n+1} \tag{6.19}$$

so the step mentioned above concerning evaluating the asymptotics of a certain y-integral is redundant. In the n-range between 4 and 100 we have seen numerically that the minimal beta-sequence appears to be relatively smooth, apart from a small wiggle near t-jump points, and to be well represented by a simple logarithm with coefficient one

$$b_n \sim \log(n/n_0) \tag{6.20}$$

There are t-jumps at the values n = 6, 17, 40 and 87, as well as further on, but overall it appears that the $t_n$ sequence also grows monotonously and logarithmically. In fact, the ratio $u_n = t_n/b_n$ is seen to vary in the n-interval 10 to 100 between 5.96 and 4.54, and although the $u_n$ sequence also features jumps, the data strongly suggest that it may approach a constant $u_n \to u$ for large n. These features are of course only indicative but they are of crucial importance for any attempt to solve the asymptotic difference equations for $b_n$ and $t_n$. In fact, instead of trying to solve the increment equations by some brute force method, which appears to be rather complicated because of the highly non-linear structure, it is preferable to approach them by trying an Ansatz and then to narrow down different choices, which subsequently can be shown to be fully consistent. If this can be done, it must indeed be "the" solution since we have seen that there exists a unique minimal beta-sequence.

Let us now focus on the difference equation for $b_n$ written in the form

$$\Delta b_n = -B_{n+1}(b_n) \frac{P^+_{n+1}(t_n/b_n)}{b_n \tfrac{\partial}{\partial b} \Xi_{n+1}(t_n, b_n)} \tag{6.21}$$

where

$$B_n(b) = \frac{4}{n+1}(1-e^{-b})^{n+1} \tag{6.22}$$

The remarkable feature of this beta-increment equation, when inserting the exact numerical $b_n$ and $t_n$ values, is that it incorporates a positive factor $B_{n+1}$ which decreases at least like 1/n and a ratio of oscillatory terms which have opposite signs, except possibly around the t-jump points where they are anyway both close to zero [10]. Therefore the beta-increment is mostly positive as expected since $b_n$ should globally be an increasing sequence.

We recall that the denominator term with the beta-derivative can be decomposed in two terms, the first being a sum of Pochhammer polynomials times a $b_{n+1}(b_n)$ factor which is very small because being O(1/n). The second term which is a finite sum of Pochhammer polynomials with $B_k^{'}$ coefficients always has the opposite sign of the numerator and overall dominates the first term, except perhaps near the t-jump points when this sum is close to zero as well. This second term was argued earlier to be close to $\Xi^*(t_n)$ but the numerical analysis shows that this is in fact the case for a t-value shifted by a small jump $t_n$, which adopting a physical oscillator terminology might more figuratively be called a phase shift, and which for large n appears to approach a constant $t$ ~ 1.12, the asymptotic size of the t-jump. We thus observe asymptotically that actually

$$b_n \tfrac{\partial}{\partial b} \Xi_{n+1}(t_n, b_n) \approx \Xi^*(t_n + t) = 4\cos(t_n + t) \qquad (6.23)$$

when evaluating $\Xi^*(t) = \Xi(t) + t\,\Xi^{'}(t)$ with $\Xi(t) = 4\sin(t)/t$. The logic of a simple cosine function appearing at this point can also be seen to be related to the polynomial Euler formula discussed in the conclusions section, and as well to the fact that the $w$-derivative of $\sin(wt)/t$ is precisely $\cos(wt)$.

We now see more clearly how the different pieces of the asymptotic analysis start to fall in place. The numerator term in the beta-increment equation is asymptotically a Pochhammer polynomial of the oscillatory trigonometric form $\cos(t_n/b_n \log(n) + j(t_n/b_n))$ which is zero at the t-jump points n=n(k) enumerated by

$$t_n / b_n \log(n) + j(t_n / b_n) = (1/2 + k)p, k = 1, 2, .. \qquad (6.24)$$

The denominator term is oscillatory of the form $\cos(t_n + t_n)$ and has the opposite sign, also being zero at the jump points, so

$$t_n + t_n = (1/2 + p)k, k = 1, 2, .. \qquad (6.25)$$

In fact, asymptotically the trigonometric phases should be equal modulo an odd multiple of $p$

$$t_n / b_n \log(n) + j(t_n / b_n) = t_n + t_n + p \qquad (6.26)$$

and since $u_n$ and $t_n$ are bounded, it follows from this that for large n we must have $t_n/b_n \log(n)$ ~ $t_n$ and therefore that $b_n$ ~ $\log(n)$ with the proportionality factor being exactly one, in accordance with the anticipated form $b_n = \log(n/n_0)$. If our full Ansatz, including $u_n \to u$ and $t_n \to t$, is correct then the equation

$$u \log(n_0) + j(u) = t + t \qquad (6.27)$$

relating the constant terms in the trigonometric phases must also be fulfilled.

But if $b_n = \log(n/n_0)$ asymptotically, then it follows that the beta-increment must satisfy the simple equation

$$\Delta b_n = \frac{1}{n} \tag{6.28}$$

which implies that all the factors on the right hand side of the general beta-increment equation must evaluate to exactly 1/n (or equivalently 1/(n+1) in the asymptotic limit). This provides a non-trivial constraint on the Ansatz which allows to determine the so far unspecified constant u. Since the trigonometric numerator and denominator terms cancel exactly because of the phase equality, we can collect the remaining (constant) factors from the expressions for $P_{n+1}^+$ and $B_{n+1}$ to get the overall factor

$$\frac{(1-n_0/n)^{n+1}}{r(u)} \approx \frac{\exp(-n_0)}{r(u)} \tag{6.29}$$

which must be exactly equal to one.

This, together with the above constant phase equation and the defining equations for r(u) and $j(u)$, leads to the following transcendental equation for u:

$$u \log(\log(\frac{1}{|\Gamma(1-iu)|})) + Arg(\Gamma(1-iu)) = t + p \tag{6.30}$$

where $t \sim 1.12$ as stated earlier. This equation has the solution $u \sim 4.4$ and from this we can then deduce that $n_0 \sim 5.10$ and $\log(n_0) \sim 1.63$ so that asymptotically $b_n \sim \log(n) -1.63$. This appears to be in excellent accordance with the numerical fit presented earlier for the case at hand. One may remark here that the argument could be modified slightly so as to determine $n_0$ from a numerical initial condition for the $b_n$ sequence which would then instead allow for evaluating the constant $t$.

We also note that the calculated asymptotic value $u \sim 4.4$ is in good agreement with the numerical data showing that u already approaches 4.5 in the n range up to 100. Moreover, we are now able to find asymptotically the values of n for which the t-roots jump:

$$n(k) = n_0 \exp(-(p+t)/u) \exp((1/2+k)p/u) \tag{6.31}$$

In particular we can extract the n-factor between successive jumps $\exp(\pi/u) \sim 2.1$ which is seen to already agree reasonably well with the low n jump values found for n = 6, 17, 40 and 87.

We have thus shown how to find a consistent formal solution to the asymptotic beta-increment equation which also appears to agree quite well with the available numerical data. Let us now discuss some more topological aspects of the discussed mechanism involving double roots of the polynomial approximants. For the lowest possible n=4 the polynomial approximant has a double root at $t_4 = 3.46$ for $b_4 = 0.074$ representing a local minimum, i.e. the second derivative is positive. Moving to n=5 the extra polynomial term has the effect of moving the graph upwards in the first quadrant because $P_5^+(u_4) > 0$. The beta-increment equation has the effect of increasing $b$ slightly so as to move closer to have a zero of $\Xi_5$, while the t-increment equation has the effect of moving t closer to a zero of $\Xi_5'$. But when we go to n=6, we see that because of $P_6^+(u_5) < 0$, the graph now moves below the axis, instead of above the axis as was previously the case, which signals that a t-jump is due and that the local minimum will become a local maximum.

The t-increment equation has the problem generally associated with the local Newton method, that the larger scale topology of the problem may not always be captured correctly. The numerical approach to finding the minimal $b$ and t-sequences therefore has to incorporate a wider extremum search routine. Alternatively it suffices to note that a change of sign of $P_{n+1}^+(u_n)$ provides a signal for t to jump by a small amount to be close to the next higher extremum of the polynomial approximant, which also happens to be close to the next zero of $\Xi^*(t)$ as seen above. We may also look at the problem of finding the minimal beta-sequence in terms of disentangling the changing topology and extremal points of a sequence of two-dimensional polynomial surfaces defined by $\Xi_n(t, b)$ for n = 4, 5…, incorporating the basic feature of having only distinct real roots in t for β = const sections when $b$ is sufficiently large.

As a final remark about the sin(t)/t case, let us restate that on the large n scale we have $t_n = u b_n$ and the t-jumps are asymptotically of constant size $t$. The t-increment equation can be used to evaluate the t-sequence for n-values between the jump values, but we will not here go into details about how to solve this equation which is a little more complex than the beta-increment equation. The t-increments are always positive because of the opposite signs of the Pochhammer polynomial term in the numerator and the second derivative term in the denominator, which indicates whether we are at a local minimum or local maximum of the polynomial. It also happens that the t-derivative of the $\log(Q_{n+1})$ factor is always positive.

## 7. Asymptotic analysis of the exponential reference cases

The full complexity of disentangling the asymptotic properties of the minimal beta-sequence is encountered when $A_I(x)$ has a genuine exponential decrease and the corresponding $\Xi(t)$ has infinitely many zeros in t, like for the reference Bessel function case. The crucial difference with the sin(t)/t case is threefold: Firstly, the number of positive $\Xi$-zeros below T has a growth rate faster than linear. Secondly, the b-coefficients of the polynomial expansion cannot be evaluated explicitly, but require asymptotic expansion of a non-trivial integral. Thirdly, since the convergence rate of the polynomial approximants is sub-logarithmic, the interest is focused on whether the growth rate of the minimal beta-sequence could possibly be sub-logarithmic as well.

Let us recall that for $\Xi(t) = 2K_{it/2}(2)$ we can easily find the low n values of the minimal $b$ and t-sequences, say starting with $b_{10} \sim 2.12$ and $t_{10} \sim 7.07$, by moving to higher and higher n using the $b$ and t-increment equations. Thus we can find $b_{100} \sim 5.57$ and $t_{100} \sim 8.11$ and note that the $b_n$ sequence has a smooth sub-logarithmic behaviour, while the $t_n$ sequence appears to move extremely slowly towards the lowest zero of $\Xi(t)$ located at t ~ 8.85. Since this simple state of affairs continues at least as far as n = 2 million [10], it is warranted to investigate a different case which might have a better chance of showing some of the peculiar minimal sequence features already found in the sin(t)/t case.

More generally, if $A_1(x) = \exp(-a(x+1/x))$ then we have $\Xi(t) = 2K_{it/2}(2a)$. Instead of the above case of a=1, we will now examine the case of choosing the much smaller value a=0.005 and immediately we note that the situation changes radically. Starting at n=4 we find $b_4 \sim 0.044$ and $t_4 \sim 1.34$ which is already located between the first and second zero of this $\Xi(t)$, at t ~ 1.29 and 2.47, respectively. Then following some small increases of $b$ and t, at n=7 we see t jump from 1.43 to 1.92 at n=8. Subsequently there are t-jumps at n = 32 and 161, as well as further on (see [10] for further details).

In this case of a=0.005, we recover all the peculiar sequence features of the previous section. Examining again the right hand side of the beta-increment equation for the exact $b_n$ and $t_n$ values, one notes that the Pochhammer term in the numerator and the beta-derivative term in the denominator both behave in an oscillatory manner, having mostly opposite signs, except eventually near the t-jump points where these terms are close to zero. Moreover it is seen that at the jump point, then t jumps from a value $t_n$ by a quantity $t_n$ to $t_n + t_n$ which is close to a zero of the function $\Xi^*(t)$. This function has a similar oscillatory behaviour as the beta-derivative term modulo the phase shift $t_n$.

The question arises immediately about why we see such an interesting minimal beta-sequence structure at quite low n for the case of a=0.005 while we do not see anything of the kind for the case a=1 as far out as n being of the order one million. The answer is simply that for the exponential cases of the A(x) family members, then are delicate logarithmic dependencies of the parameters, like the a-parameter introduced above, and that some features may be more or less visible depending on the actual numerical size of the parameters. Consequently a purely numerical analysis coupled with our finite computing power would logically be expected to provide a quite partial understanding of the true asymptotics, which only an appropriate formal asymptotic analysis might overcome.

Thus it is warranted to try to disentangle the asymptotic behaviour of the minimal beta-sequence by investigating the formal mechanism exposed in the previous section related to the asymptotic properties of all the terms entering into the right hand side of the beta-increment equation. The key points to understand are the asymptotic properties of the $B_{n+1}(b_n)$ coefficients and their expression in terms of a non-trivial integral, as well as the asymptotic properties of the beta-derivative term in the denominator, which amounts to analyzing the asymptotic form of the function $\Xi^*(t)$. The latter one is quite well-known so let us just state the result that for $\Xi(t) = 2K_{it/2}(2a)$, we have the large t asymptotic form

$$\Xi(t) \approx \frac{4}{\sqrt{pt}} \exp(-pt/4) \cos(t/2 \log(t/(2ae)) - p/4) \tag{7.1}$$

The zeros of this function are given explicitly by

$$t_k = 2ae \text{ sublog}\left[\exp\left[(k-1/4)p/(ae)\right]\right], k=1,2,... \tag{7.2}$$

which are quite close to the zeros of the exact $\Xi(t)$ even for small k. The number of positive zeros in the interval [0,T] is therefore

$$N(T) \sim T/(2p)\,\log(T/(2ae)) \tag{7.3}$$

whereas in the simple periodic sin(t)/t case this number of course only grows linearly.

For large t we thus find

$$\Xi^*(t) \approx -2\sqrt{\frac{t}{p}}\exp(-pt/4)\log(t/(2ae))\sin(t/2\log(t/(2ae)) - p/4) \tag{7.4}$$

It remains to evaluate the large n behaviour of the $B_{n+1}(b_n)$ coefficient

$$B_{n+1}(b_n) = 4\int_0^1 dy\,\exp(-a(y^{-e(n)} + y^{e(n)}))(1-y)^{n+1} \tag{7.5}$$

where as earlier $e(n) = 2/b_n$ for short. We immediately note that the second term in exponent of the integrand can be dropped because it gives a subdominant contribution in the integration range near y=0 which is responsible for defining the large n behaviour of $B_{n+1}$. This is in stark contrast to the importance of this exponential term in assuring the x-inversion invariance of the corresponding $A_1(x)$, which in turn is instrumental for $\Xi(t)$ having an infinite number of zeros and decreasing exponentially for large t. This point will also be seen to be quite relevant when discussing the incomplete gamma function cases below.

The Laplace method can be applied directly to evaluate this integral and we will therefore have to study the extremal point of a certain function

$$f(y) = (n+1)\log(1-y) - ay^{-e(n)} \tag{7.6}$$

The precise asymptotic form of the $B_{n+1}$ integral depends on the decrease rate of $e(n)$ (or equivalently the growth rate of $b_n$) but it will always contain a factor of the of the form $\exp(-a\,n^{e(n)})$. If we suppose that ε(n) decreases as slowly as $e(n) \sim (\log(n))^{-d}$ with 0<$d$ <1, then this exponential factor will go very fast to zero for $n \to \infty$, and in fact so fast that it will be seen to dominate all other factors in the right hand side of the beta-increment equation. Therefore there is little chance that the corresponding $b_n$ would constitute a good Ansatz for a self-consistent solution to this equation, so we will leave aside this possibility for now and comment on it again below.

On the other hand, if $e(n)$ decreases as fast as $e(n) \sim 1/\log(n)$, then the above exponential factor would approach a constant and would therefore play an insignificant role in balancing the different factors in the increment equation. For ε(n) = 1/sublog(n) then $\exp(-a\,n^{e(n)}) = \exp(-a/e(n))$ which is not decreasing too fast for large n. We therefore argue that the interesting range for a suitable Ansatz for $b_n$ in the beta-increment equation is 1/sublog(n)≥ $e(n)$ >>1/log(n). If $e(n)$ fulfils this constraint then the asymptotic form of $B_{n+1}$ can be evaluated by the Laplace method to be

$$B_{n+1}(b_n) \approx \frac{1}{n}\sqrt{2\pi a e(n) n^{e(n)}} \exp(-a n^{e(n)} - a e(n) n^{e(n)}) \qquad (7.7)$$

an expression in which all n-dependent terms, as well as the overall constant, will be seen below to acquire a special importance, as in the final analysis of the beta-increment equation for the sin(t)/t case.

Just like in the sin(t)/t case we observe that the beta-derivative term of denominator in the right hand side of the beta-increment equation is simply $\Xi^*(t_n + t_n)$, where $t_n$ is related to the size of the t-jump. In the sin(t)/t case the density of $\Xi(t)$ zeros is constant which is why $t_n$ was seen to approach asymptotically a non-zero constant $t$. However in the exponential case at hand, the density of zeros grows logarithmically and therefore the size of the t-jumps will go to zero, so that $t_n$ will actually here tend to $t = 0$.

The requirement of having a self-consistent solution of the beta-increment equation, as well as the constraints indicated by the stated features of the numerical data, implies first of all that asymptotically the trigonometric term of the numerator and denominator, i.e. $\cos(t_n/b_n \log(n) + \pi (t_n/b_n))$ and $\cos(t_n/2 \log(t_n/(2ae)) + \pi/4)$, cancel each other. A priori, it is a open question whether $u_n = t_n/b_n$ approaches a non-zero constant for large n, like in the sin(t)/t case, or tends towards zero, as the numerical data seems to suggest. In any case, the equality of the trigonometric phases (modulo an odd multiple of $\pi$) clearly implies that the leading terms must match, that is $t_n/b_n \log(n) \sim t_n/2 \log(t_n/(2ae))$ so that we get the simple asymptotic relation between the beta and t-sequences

$$\log(n) = b_n/2 \, \log(t_n/(2ae)) \qquad (7.8)$$

Using the definition $e(n) = 2/b_n$ this can be recast in the form

$$n^{e(n)} = t_n/(2ae) \qquad (7.9)$$

which when inserted in the beta-increment equation makes the asymptotic balancing of the different terms quite transparent, as will be seen below. Let us keep in mind that we are now trying to find a suitable Ansatz for $b_n$ fulfilling the constraint $1/\text{sublog}(n) \geq e(n) \gg 1/\log(n)$, which conveniently may be parametrized by the implicit equation

$$\log(n) = 1/e(n) \, g(1/e(n)) \qquad (7.10)$$

where g is a slowly increasing function. If $g(1/e(n)) = \log(1/e(n))$ then of course $1/e(n) = \text{sublog}(n)$, but other cases where g increases even slower that the logarithmic function could potentially be of interest. For example, if $g(1/e(n)) = \log(\log(1/e(n)))$ then, depending on the introduction of additional constant factors or sub-dominant terms in this definition, then we get a type of extra-large sub-logarithmic growth of $1/e(n)$ which we will here generically denote by the function name sublogxl(n), suggested by the relation sublog(n) << sublogxl(n) << log(n).

If $1/e(n) = \text{sublog}(n)$ then the asymptotic beta-increment equation to be satisfied is

$$\Delta b_n = \frac{2}{n(\log(b_n/2)+1)} \approx \frac{2}{n\log(b_n)} \qquad (7.11)$$

and in the general g case we must have asymptotically that

$$\Delta b_n = \frac{2}{ng(b_n)} \qquad (7.12)$$

Following these remarks, we now know that we are looking for a $b_n$ Ansatz which should be such that the right hand side of the beta-increment equation contains an overall 1/n factor, as well as an additional slowly decreasing 1/g-factor. In order to see how the different terms may match asymptotically, we first observe that the exponential terms in the numerator and the denominator are the most dominating terms, and the Ansatz should first of all be such that these terms cancel. On one hand, the numerator $B_{n+1}$ term is seen to contain the dominant exponential factor $\exp(-a\, n^{e(n)}) = \exp(-t_n/(2e))$, when using the above relation between $b_n$ and $t_n$. On the other hand, the denominator term contains the exponential factor $\exp(-pt_n/4)$, which decreases faster than the numerator exponential, so they cannot possibly cancel unless some other factor would be of a comparable asymptotic size.

Clearly the other $t_n$ factors in $\Xi^*(t_n)$ are not relevant in this preliminary context. Only the $(e(n))^{1/2}$ factor in the $B_{n+1}$ term might decrease sufficiently fast so as to be comparable to an $\exp(-t_n)$ type of factor. But this means that asymptotically we must have that $\log(b_n)=O(t_n)$, which implies that $t_n$ should grow much slower than $b_n$, in accordance with the numerical data suggesting that $u_n = t_n/b_n$ tends towards zero. This preliminary analysis leads us to propose the following Ansatz for the beta-increment equation:

$$\log(n) = b_n/2\, (\log(s/(2ae)\, \log(b_n/m)) + q/(s\, \log(b_n/m)) \text{ and } t_n = s\, \log(b_n/m) \qquad (7.13)$$

where $s, m$ and $q$ are constants to be determined. We note here that the proposed $b_n$ is of sublogxl(n) type as defined above.

At this point we may pause to remark that it can now be seen clearly that the previously considered possibility of the minimal $b_n$ sequence increasing as slowly as a power of a logarithm $b_n \sim (\log(n))^d$ with δ<1 cannot possibly respect the delicate balance of different terms in the beta-increment equation, simply because the $\exp(-a\, n^{e(n)})$ factor will dominate all other terms so that no cancellation of leading behaviours is possible. This implies that the previously stated numerical fit to the exponential $b_n$ data where $d \sim 0.6$ is unlikely to represent the true asymptotic behaviour of $b_n$, even if the low n fit was seen to be quite excellent. This type of sub-logarithmic fit was of course just the simplest proposal for a parameterization of the data, but it is in fact possible to propose an equally good alternative parameterization of the sublogxl(n) type.

To pursue the asymptotic analysis, we will show that the above Ansatz actually provides a consistent asymptotic solution of the beta-increment equation. Since $u_n = t_n / b_n$ tends to u=0 we can use r(0)=1 and $j(0) = 0$ to simplify somewhat the calculations, so let us here restate the beta-increment equation as

$$\Delta b_n = -B_{n+1}(b_n) \frac{\cos(t_n / b_n \log(n))}{\Xi^*(t_n)} \tag{7.14}$$

The trigonometric phases of the numerator and denominator terms of the right hand side of this equation must be equal modulo an odd multiple of π, which can be found from the initial condition of the beta and t-sequences, so more precisely the following equation has to be fulfilled

$$t_n / b_n \log(n) = t_n / 2 \log(t_n / (2ae)) + 21p/4 \tag{7.15}$$

This implies a small correction to the relation stated earlier

$$n^{e(n)} = t_n /(2ae) \exp(21 p /(2 t_n)) \sim t_n /(2ae) + 21 p /(4ae) \tag{7.16}$$

and it also follows that in the Ansatz above that $q = 21p/4$. Combining the leading exponential factors in the beta-increment equation we get

$$\exp(-an^{e(n)}) \exp(pt_n/4)\sqrt{e(n)} \approx \sqrt{\frac{2}{m}} \exp(-q/e - t_n/(2e) + pt_n/4 - t_n/(2s)) \tag{7.17}$$

If we now impose

$$s = \frac{1/2}{p/4 - 1/(2e)} \approx 0.831 \tag{7.18}$$

then these 3 factors altogether just leave the constant factor

$$\sqrt{\frac{2}{m}} \exp(-q/e) \tag{7.19}$$

Next we note that the two factors which behave like a power of $t_n$ cancel each other automatically:

$$(n^{e(n)})^{1/2} (t_n)^{-1/2} \approx (2ae)^{-1/2} \tag{7.20}$$

By a some kind of miracle, the only n-dependent factor which remains asymptotically is the logarithmic factor in the denominator term $\Xi^*(t_n)$:

$$\log(t_n /(2ae)) \approx \log(\log(b_n)) \tag{7.21}$$

which is precisely the g-function type of factor needed for self-consistency of the beta-increment equation when $b_n \sim $ sublogxl(n) as proposed in the Ansatz! It only remains to collect all constant factors multiplying this term and the 1/n term, and verify that they can actually amount to exactly 2. This leads to an additional final equation which determines the remaining parameter $m$ of the Ansatz:

$$m = \frac{2p^2}{e} \exp(-2q/e) \approx 3.89 \cdot 10^{-5} \qquad (7.22)$$

In summary, we have seen that beta-increment equation can be fulfilled in the general case of $\Xi(t)$ being of Bessel function form by a minimal beta sequence $b_n$ of sublogxl(n) type with $t_n \sim \log(b_n)$. As in the sin(t)/t this is a large scale asymptotic statement which includes the t-jumps. The more precise behaviour of $t_n$ in between the jumps can be deduced from the t-increment equation, but we will not here go into detail about this technical point.

Let us now discuss briefly the other exponential reference cases concerning $\Xi(t)$ being a symmetrized incomplete gamma function. One of the good comparational reasons for specifically considering this case is, that if one retains only a finite number of terms in the elliptic theta function expression for the Riemann A(x), then the corresponding $\Xi(t)$ will precisely be a finite sum of incomplete gamma functions.

We may recall that for $A_I(x) = \exp(-x)$ then we obtain the symmetrized incomplete gamma function $\Xi(t) = \Gamma(it/2, 1) + \Gamma(-it/2, 1)$. As in the Bessel case it is easy to find the low n values of the minimal $b$ and t-sequences, say starting with $b_{10} \sim 2.76$ and $t_{10} \sim 8.32$ and moving further to note that the $b_n$ sequence has a smooth sub-logarithmic behaviour, while the $t_n$ sequence appears to increase much more slowly. This situation continues at least as far as n=2 million [10], and again one may ask the question whether some other parameter choice might possibly show some of the peculiar minimal sequence features found in the sin(t)/t case.

A key difference with the Bessel case is of course that the $\Xi(t)$ is positive with no (real) zeros, except the one at infinity. This does of course not prevent the polynomial approximants from having roots but, if the approximants converge, one would simply expect all roots to move to infinity with increasing n. Another importance difference is that the Bessel $\Xi(t)$ decreases exponentially in t while the symmetrized incomplete gamma function $\Xi(t)$ decreases only quadratically for large t:

$$\Xi(t) \approx \frac{8}{et^2} \qquad (7.23)$$

Let us now consider the more general case $A_I(x) = \exp(-ax)$ with a>0 where we also obtain a symmetrized incomplete gamma function

$$\Xi(t) = a^{-it/2} \Gamma(it/2, a) + a^{it/2} \Gamma(-it/2, a) \qquad (7.24)$$

We can find the asymptotic behaviour of this $\Xi(t)$ by using the general relation

$$\Gamma(z,a) = \Gamma(z) - a^z e^{-a} \sum_{k=0}^{\infty} a^k \frac{\Gamma(z)}{\Gamma(z+k+1)} \Gamma(z,a) \tag{7.25}$$

where the first term asymptotically gives rise exactly to the exponentially decreasing expression of the Bessel case, whereas the leading terms of the sum give rise to a quadratically decreasing term, so that altogether

$$\Xi(t) \approx \frac{8a}{e^a t^2} \tag{7.26}$$

If a is sufficiently small, actually if a is smaller than 0.23, then $\Xi(t)$ will start having zeros, and it is therefore interesting, like in the Bessel case, to see what happens if a is quite small, say a=0.01. In this case there is exactly 8 zeros located at 1.49, 2.81, 4.04,…Again we note that the beta-sequence structure changes radically, because starting at n=4 we find $b_4 \sim 0.054$ and $t_4 \sim 1.53$ which is already located between the first and second zero of $\Xi(t)$. Then following some small increases of $b$ and t, at n=7 we see t jump from 1.62 to 2.12 at n=8. Subsequently there are t-jumps at n = 36 and 248, as well as further on.

The interesting question is now how to understand the behaviour of the beta and t-sequences for low n, as well as their ultimate asymptotic behaviour, in terms of the beta and t-increment equations. For low n we see numerically exactly the same features of the terms of these equations as observed in the sin(t)/t and Bessel case. That is, the $P_{n+1}^+(u_n)$ term in the numerator and the $\Xi^*(t_n)$ term in the denominator have opposite oscillatory behaviour. But this situation can only continue, as $t_n$ increases, in a finite range of n values, because $\Xi^*(t)$ just has a finite number of zeros. From a certain point on, then $\Xi^*(t)$ will not have more zeros, and so $\Xi^*(t) \sim -1/t^2$.

We thus see that the asymptotic structure of the problem is somewhat different in this case, but nevertheless a number of features remain similar to the Bessel case. In the beta-increment equation, for $A_I(x) = \exp(-ax)$ the $B_{n+1}$ coefficient still has the same asymptotic form and if we make the same Ansatz as before that $u_n$ tends to zero, then the beta-increment equation will asymptotically have the form

$$\Delta b_n = B_{n+1}(b_n \frac{e^a t_n^2}{2a} \cos(t_n / b_n \log(n)) \tag{7.27}$$

The t-increment equation can now be simplified asymptotically because the logarithmic derivative of the beta-derivative term will be sub-dominant and we therefore have

$$\Delta t_n = -\frac{b_{n+1}(b_n) P_{n+1}^{+'}(t_n / b_n)}{\Xi"_{n+1}(t_n, b_n)} \approx B_{n+1}(b_n) \frac{3e^a t_n^4 \log(n)}{a b_n^2} \sin(t_n / b_n \log(n)) \tag{7.28}$$

For low n, the t-jump mechanism, as in the Bessel case, was related to having synchronized numerator and denominator terms, but for the large n asymptotic limit the jump mechanism is now somewhat different. A t-jump is still signaled in the beta-increment equation by the cosine argument coming close being $p/2$ modulo $p$, but since both the beta and the t-increments are globally positive, then the t-jumps must be tuned to assure that the trigonometric phase $u_n \log(n)$ is confined to always being in the first quadrant.

After this argument about how the t-jump mechanism works in this case, one may try to understand how the asymptotic behaviour of the different factors may cooperate to produce a self-consistent solution to the increment equations. The main point of this analysis is, as before, that the $b_n$ growth rate should be sufficiently fast so to as assure that the factor $\exp(-a\, n^{e(n)})$ in the $B_{n+1}$ term is of comparable asymptotic size as the other factors. This imposes again the requirement that $b_n$ should be of sublogxl(n) type. But since the other leading factors are powers in $b_n$ and $t_n$, we are now led to an additional Ansatz of a relation between $b_n$ and $t_n$ of the form $t_n \sim (b_n)^d$ with $\delta<1$. We will not here go into further technical details of this case, since it is not central to the analysis of the Riemann case in the following section. In summary, we have argued that the $b_n$ growth rate is of sublogxl(n) type like the Bessel case, and also that there are qualitative differences which are responsible for the available numerical data indicating a slightly different type of sub-logarithmic growth.

## 8. Asymptotic analysis of the Riemann Xi-function case

The most well-known bona fide approximation to the Riemann $\Xi$ function [24] arises when keeping only the leading term in the elliptic theta function representation of the Riemann A(x), and then modifying the exponent by substituting x+1/x instead of x, i.e.

$$A_I(x) = 2p^2(x+1/x)^{9/4} \exp(-p(x+1/x)) \tag{8.1}$$

This choice leads to a Bessel function expression of the form

$$\Xi(t) = p^2 (K_{9/4+it/2}(2p) + K_{9/4-it/2}(2p)) \tag{8.2}$$

which asymptotically evaluates to

$$\Xi(t) \approx 2^{-5/4} p^{1/4} t^{7/4} \exp(-pt/4) \cos(t/2 \log(t/(2pe)) + 7p/8) \tag{8.3}$$

We can apply the beta-increment analysis to this case almost exactly as it was done earlier for the Bessel reference case, with the special choice of $a = p$. The only essential difference in the asymptotic form of the $\Xi(t)$ is the constant $7p/8$ versus $-p/4$ in the cosine argument and the factor $t^{7/4}$ versus $t^{-1/2}$. This difference is entirely due to the extra factor $x^{9/4}$ in $A_I(x)$ as compared to the reference case. It is quite easy to see that if such an extra x power factor is included in the integral defining $B_{n+1}(b_n)$, then the asymptotic expression for this term in the beta-increment equation also acquires an extra factor $(n^{e(n)})^{9/4} \sim (t_n)^{9/4}$, just like $\Xi(t)$ and $\Xi^*(t)$. This means that the delicate asymptotic cancellation mechanism between numerator and denominator terms is still operational with the same kind of Ansatz. Consequently, in this more complex $A_I(x)$ case, we also have that $b_n \sim 2$ sublogxl(n) is a consistent solution to the beta-increment equation.

This bona fide approximation to the Riemann case has the property of the number of positive real zeros in the interval [0,T] being

$$N(T) \sim T/(2p) \log(T/(2pe)) \tag{8.4}$$

This is identical to one of the consequences of the Riemann Hypothesis. But the above asymptotic $\Xi(t)$ behaviour is not the same as for the Riemann $\Xi(t)$ for the simple reason that the approximation neglects terms in the sum expression for the Riemann A(x). However, if we add the neglected exponential terms, still keeping the leading $x^{9/4}$ factor in $A_1(x)$ and making the x+1/x substitution, then we get the much better approximation

$$\Xi(t) = p^2 \sum_{k=1}^{\infty} k^4 (K_{9/4+it/2}(2pk^2) + K_{9/4-it/2}(2pk^2)) \tag{8.5}$$

The asymptotic expansion of this expression for large t is identical to the asymptotic form given above except that the sum, before taking the real part of the expression, produces an extra factor

$$\sum_{k=1}^{\infty} \frac{1}{\sqrt{k}} k^{-it} \tag{8.6}$$

But this is just the formal expression for $V(1/2+it)$, so the considered extension of the bona fide approximation to the Riemann $\Xi(t)$ evaluates in fact asymptotically to precisely what it should be! There is of course a direct way of finding the asymptotics of the exact Riemann $\Xi(t)$ using the relation

$$\Xi(t) = -1/2(t^2 + 1/4)p^{-1/4-it/2}\Gamma(1/4+it/2)z(1/2+it) \tag{8.7}$$

The asymptotic form of this expression is

$$\Xi(t) \approx (p/2)^{1/4} t^{7/4} \exp(-pt/4) \exp(it/2 \log(t/(2pe)) + i7p/8) z(1/2+it) \tag{8.8}$$

and with the standard introduction of the real quantity Z(t), having the opposite sign of $\Xi(t)$, and the phase of $V(1/2+it)$ defined through the relation

$$Z(t) = e^{iJ} z(1/2+it) \tag{8.9}$$

then we can rewrite this as

$$\Xi(t) \approx -(p/2)^{1/4} t^{7/4} \exp(-pt/4) Z(t) \tag{8.10}$$

As in the Bessel reference case, Z(t) has an infinite number of zeros, but a significant difference is that $|Z(t)| = |V(1/2+it)|$ is unbounded in contrast to being a simple cosine. Nevertheless we can apply the beta and t-increment equations in the same way as in the reference case. We recall that it was found that $b_{10} \sim 3.23$ and $t_{10} \sim 10.63$, and continuing to higher n, we find a very smooth sub-logarithmic behaviour of $b_n$ and an even slower smooth behaviour of $t_n$. Moreover, this smooth behaviour persists

at least up to n=2 million [10]! However, we know from the Bessel case, thanks to the fact that it was possible to examine different a-parameter cases, that for higher n then at some point we must eventually start to see t-jumps, and that this process will continue ad infinitum.

The central question is now whether the mechanism involving the beta and t-increment equations, which in the Bessel case was seen to lead to $b_n \sim 2$ sublogxl(n) being a consistent solution, can be applied with a similar result in spite of Z(t) having a complex behaviour which is quite far from being completely understood. This element of uncertainty suggests that the analysis of the minimal beta-sequence for the Riemann case might also be a rather subtle matter. However, the most significant element in this analysis is still the very fact that there exists a unique minimal beta-sequence solution to the increment equations, and that the delicate asymptotic balancing of different terms in these equations implies some strong constraints.

The very close analogy between the Bessel case and the Riemann case in what concerns the leading asymptotic factors in the numerator and denominator terms, as well in the synchronization mechanism between the oscillatory numerator and denominator terms, forces the search for a consistent solution towards the same sublogxl(n) type of growth for the minimal beta sequence. This circumstance imposes strong constraints on the behaviour of the rest of terms, including on the asymptotic behaviour of Z(t). In other words, if the quite elusive Z(t) function would not have a certain property, which might be a hitherto unproven property, then there would not exist a consistent solution to the increment equations. Consequently, due to the fact that a unique solution of these equations does exist, then we would be able to infer that some specific property of Z(t) must indeed be true.

Let us now consider a first significant example of this type of consistency reasoning which is connected with the asymptotic distribution of the zeros of Z(t). The first key observation in this respect by Selberg [23] was that the number of real zeros $N_0(T)$ of Z(t) grows faster that linearly like

$$N_0(T) \sim d\ T/(2p)\ \log(T) \tag{8.11}$$

where $d \leq 1$ is a very small number. This implies that a finite proportion $d$ of all zeros N(T) of Z(t) (i.e. including possible complex ones) lie on the critical line. Subsequently (see [6, 12]}, this result has very laboriously been improved using the so-called "mollifier" technique, so that we now can state, using a common jargon, that about 41% of all zeros lie on the critical line. These zeros can also be shown to be simple. If the Riemann Hypothesis is correct, then of course we have the implication that 100% of all zeros lie on the critical line.

However, it turns out that this latter conclusion appears quite naturally in the framework of the minimal beta-sequence analysis using the following argument. Like in the Bessel case we should be looking for beta and t-sequence solutions for which $u_n \to 0$, so the Pochhammer polynomial in the numerator of the right hand side of the beta-increment equation becomes simply $\cos(t_n / b_n \log(n))$. Moreover, the zeros of this cosine are asymptotically synchronized with the zeros of the denominator $\Xi^*(t_n)$, which has the oscillatory factor $Z'(t)$. Therefore the key asymptotic relation between $b_n$ and $t_n$, similar to the one found in the Bessel case, now reads $t_n / b_n \log(n) = d\ t_n / 2\ \log(t_n)$ which, with $e(n) = 2/b_n$ as earlier, implies

$$n^{e(n)} = (t_n)^d \tag{8.12}$$

But if there would be a strict inequality $d<1$ then immediately we run into a problem with the delicate asymptotic balance of the leading numerator and denominator terms in the beta-increment equation which applies fully to the Riemann case. Simply, the numerator exponential factor $\exp(-p\ n^{e(n)}) \sim \exp(-p(t_n)^d)$ from the $B_{n+1}$ coefficient now has a much more different asymptotic behaviour than the $\Xi^*(t_n)$ denominator factor $\exp(-pt_n/4)$, making mutual cancellation, and even "assisted" cancellation, impossible. In the Bessel case, the $B_{n+1}$ factor $(e(n))^{1/2}$ could be brought into the picture to assure the asymptotic cancellation of the most leading terms, but this would not be possible now, in contrast to the Bessel case where we know for sure that $d=1$.

Consequently, in case $d<1$ then we find a contradiction with the requirement of the existence of a minimal beta-sequence solution. We thus conclude from the above argument that $d=1$, so that 100% of all zeros must lie on the critical line. These "100%" of course do not imply that the Riemann Hypothesis is true since there might still theoretically be a few isolated zeros off the critical line. To address this question, as it will be done in the final section, a different argument based on the Hurwitz theorem of complex analysis is needed.

The above argument can be made a little more precise if we include an additional constant in the $N_0(T)$ expression:

$$N_0(T) \sim T/(2p) \log(T/T_0) \tag{8.13}$$

Similarly to the Bessel case where $T_0 = 2pe$, we now find $n^{e(n)} = t_n/T_0$ and therefore the numerator factor $\exp(-p\ t_n/T_0)$ has to be balanced with the denominator factor $\exp(-p\ t_n/4)$. This balancing would be perfect if $T_0=4$ but otherwise it is imperative to have $T_0>4$, in which case one has to include the numerator factor $(e(n))^{1/2}$ in the balance accounting, as done in the Bessel case. In the latter case, we are of course again led to the conclusion that $b_n$ must be of sublogxl(n) type. In the hypothetical case $T_0=4$, then a consistent solution of the beta-increment equation may be found which has an even slower beta growth rate between sublog(n) and sublogxl(n), but we will not here go further into analyzing this possibility, since the subsequent arguments would be similar.

The second significant example of the reasoning related to the existence of a unique minimal beta-sequence solution concerns the asymptotic growth of Z(t) and emerges when analyzing the balancing of the next to leading asymptotic terms in the beta-increment equation, i.e. the terms behaving like a power of $t_n$. The argument here starts with the observation that if $b_n$ is of sublogxl(n) type then, just like in the Bessel case, the numerator and denominator factors of degree $t^{9/4}$ will balance out exactly. But in the Riemann case the question immediately arises about the effect of the additional $Z'(t_n)$ factor in the denominator coming from $\Xi^*(t_n)$. This factor could potentially have a power like-behaviour $t^d$ like Z(t) and, superficially at least, there would seem to be nothing left to cancel it.

In this respect, let us first recall that the Lindelöf Hypothesis [24] states that for any $d>0$ then

$$Z(t) = O(t^d) \tag{8.14}$$

but this property has so far only be proved for $d$ as small as about 1/6. It is of course well-known that the Lindelöf Hypothesis follows from the Riemann Hypothesis. We may also recall that in the Bessel case, where Z(t) is essentially replaced by $\cos(t/2 \log(t/(2pe)))$, then $\Xi'(t)$ acquires asymptotically an extra factor of $1/2 \log(t)$ compared to $\Xi(t)$. In the balancing mechanism, this provides the only surviving, slowly growing factor $\log(t_n) \sim \log(\log(b_n))$ which was precisely seen to be what was needed for having a self-consistent solution to the beta-increment equation of type sublogxl(n).

Thus we observe that the requirement for a consistent solution to the beta increment calls for the presence of a denominator factor growing like log(t), but that in the Riemann case we instead have a $Z'(t)$ factor, which essentially is just known to be unbounded. The solution to this problem can be found by remarking that there are non-leading corrections to the asymptotic formula for the number of zeros of the form

$$N_0(T) \sim T/(2p) \log(T/2pe) + S(T) \tag{8.15}$$

which are big enough to play a significant role in our analysis. In fact, we know that the additional S(T) term might be as important as O(log(T)) and that S(T) is an infinitely alternating function. The previous asymptotic relation between $b_n$ and $t_n$ therefore has to be modified to read

$$n^{e(n)} = t_n/(2pe) \exp(2p S(t_n)/t_n) \sim t_n/(2pe) + S(t_n)/e \tag{8.16}$$

Consequently the numerator factor $\exp(-p\, n^{e(n)})$ gives rise to the additional factor $\exp(-pS(t_n)/e)$, which is positive, and which may be both quite small and quite big. But this means that, after balancing all explicit exponential and power factors in the numerator and the denominator of the beta-increment equation, it remains only to balance the numerator factor $\exp(-pS(t_n)/e)$ with the denominator factor $Z'(t_n)$ to produce a $1/\log(t_n)$ factor, as required by the consistency requirement.

In summary, if we impose the asymptotic growth constraint

$$|Z'(t)| \sim \log(t) \exp(p|S(t)|/e) \tag{8.17}$$

then we can assert that the solution $b_n \sim 2\text{sublogxl}(n)$ is perfectly consistent and that there is no other way of assuring consistency. It should here be reemphasized that this growth rate for the minimal beta-sequence emerged as a direct consequence of the asymptotic balancing of the leading terms, independently of the presence of Z(t) and S(t) because of known bounds on the growth rate of these. Subsequently the above relation between Z(t) and S(t) imposed itself simply because there must exist a unique solution to the beta-increment equation.

The above relation is relatively weak in the sense that it does not explicitly imply specific individual properties of Z(t) and S(t). Moreover, since the derivative of the $V(1/2+it)$ phase is of the order $1/2 \log(t)$ asymptotically, it is reasonable to expect that $|Z'(t)|$ should not be very different from $\log(t)|Z(t)|$ asymptotically. Therefore the relation can be seen to be quite compatible with already established bounds for Z(t) and S(t), as well as with the ones which follow from both the Lindelöf and the Riemann Hypothesis for these quantities.

## 9. Conclusions

For the final step of our analysis we now invoke the Hurwitz theorem of complex analysis [1,18] stating, that if an analytic function is a limit of a sequence of analytic functions, uniformly convergent on a compact subset of the complex plane, then any zero of the function in the subset must be a limit of the zeros of the sequence functions. This theorem allows us to arrive at several interesting conclusions about the zeros of the entire functions $\Xi(t)$ belonging to the broad family of functions defined by an integral representation of the form

$$\Xi(t) = \int_1^\infty dx A(x) x^{-1/4} (x^{it/2} + x^{-it/2}) \qquad (9.1)$$

where A(x) is real, non-negative, positive and continuous at x=1, bounded on $[1,\infty]$ and decreasing like $\exp(-ax^b)$ with a>0 and b>0 or faster for $x \to \infty$.

It is not possible at the present stage of our understanding of this function family to transform the insights concerning the features of the polynomial approximants $\Xi_n(t, b_n)$ to $\Xi(t)$ discussed in the previous sections into a single, succinct theorem characterizing all the different family members A(x) with respect to the structure of the zeros of the associated $\Xi(t)$. However, the discussed insights do allow to draw up a quite differentiated picture of the A(x) family, in particular concerning the Riemann A(x), which amounts to a compelling confirmation of the validity of the Riemann Hypothesis.

The key insights discussed earlier, underlying our conclusions, are that any $\Xi(t)$ function can be expressed as a uniformly convergent expansion of symmetrized Pochhammer polynomials which have some unique real root properties. These properties are specifically expressed through the appearance of a real root regime for the polynomial approximants $\Xi_n(t, b)$ which allow us to infer the existence of a mapping

$$A(x) \to b_{\min,n} \qquad (9.2)$$

which associates to every individual family member A(x) a minimal beta-sequence $b_{\min,n}$ (n=1,2,..), having the property that for all $b \geq b_{\min,n}$ then the polynomial $\Xi_n(t, b)$ has only real roots. Thus the elements of the polynomial approximation sequence $\Xi_n(t, b_n)$ have real roots only for large n if $b_n$ is chosen to grow asymptotically at least as fast as $b_{\min,n}$, i.e. $b_n \geq b_{\min,n}$ for large n. In general, it is a non-trivial problem to find explicitly the minimal beta-sequence associated with a given A(x)

On the other hand, we have shown that the sequence $\Xi_n(t, b_n)$ converges to $\Xi(t)$ if $b_n$ grows slower than a certain maximal sequence $b_{\max,n}$ or possibly at the same rate, i.e. $b_n \leq b_{\max,n}$ for large n. Specifically, if $A_1(x) = 1$ with compact support then we have $b_{\max,n} = B_0 + \log(n)$ and in this case the maximum growth rate can actually be attained, but the price to be paid for this is that the convergence rate is just logarithmic. If A(x) is of genuine exponential decrease then $b_{max,n} = \log(n)$, but this growth rate cannot be attained and we must have $b_n = o(\log(n))$ for convergence.

Consequently, the Hurwitz theorem allows us to conclude that $\Xi(t)$ has only real zeros if for large n we fulfill the growth rate relation

$$b_{min,n} \leq b_{max,n} \tag{9.3}$$

for the case of A(x) having compact support, and the relation

$$b_{min,n} \ll b_{max,n} \tag{9.4}$$

for the case of A(x) being of genuine exponential decrease.

The central question is therefore what exactly is the $b_{min,n}$ growth rate for a given A(x). Numerical analysis of the minimal beta-sequences and associated t-root sequences can be performed to any desired precision and very high polynomial degree by iteration of specific difference equations for $b_n$ and $t_n$. This kind of analysis provides a first indication that the minimal growth rate for the reference case of when $\Xi(t) = 4\sin(t)/t$ is simply $\sim \log(n)$, while for the reference case of the Bessel function and the Riemann $\Xi(t)$ we find indications of a distinct sub-logarithmic behaviour.

Equally importantly, the numerical analysis provides some crucial insight into the fundamental mechanism which is behind the relatively smooth growth behaviour of the $b_n$ sequence and the more jumpy behaviour of the $t_n$ sequence. Numerical evidence can, however, only provide a partial picture of the true asymptotic behaviour of these sequences which is described by two succinct beta and t-increment equations. For the $\Xi(t) = 4\sin(t)/t$ reference case, this asymptotic analysis does confirm that the minimal beta-sequence is exactly $\log(n/n_0)$ and the Hurwitz theorem therefore allows us to conclude that this $\Xi(t)$ has only real zeros. This fact is of course already known.

For the Bessel function reference case, the asymptotic analysis shows that the growth rate of the minimal beta-sequence is of type sublogxl(n)<<log(n) and the Hurwitz theorem therefore allows us to conclude that this $\Xi(t)$ has only real zeros. This fact is of course also already known. But for the Riemann case, the asymptotic analysis actually provides information which is definitely not known, even if it may have been widely anticipated!

Firstly, the asymptotic analysis provides the information that "100%" of all Riemann $\Xi(t)$ zeros lie on the critical line, as expressed using a common jargon. Secondly, it shows that the growth rate of the minimal beta-sequence is of type sublogxl(n)<<log(n), like in the Bessel case, and it imposes a certain relation between the functions Z(t) and S(t), Consequently, we can also here apply the Hurwitz theorem to argue that the Riemann $\Xi(t)$ has only real zeros, thus providing a compelling confirmation of the validity of the Riemann Hypothesis.

The asymptotic analysis and the application of the Hurwitz theorem has been possible for some representative member of the considered A(x) family, including the most interesting member the Riemann A(x). But it has not yet been possible to disentangle the situation for all other interesting members of the family for the simple reason that it not always evident how to characterize some of the important parameters of the corresponding $\Xi(t)$ which enter into the beta and t-increment equations, and to find a consistent solution of these equations.

These important parameters include the asymptotics of $\Xi(t)$ for large real t and the distribution of its known real zeros. In a sense, we have seen that it helps if it happens that $\Xi(t)$ has an infinite number of real zeros. If there are only finitely many zeros, or no zeros, then the structure of the increment equations changes and a different mechanism for the asymptotic behaviour is at play. For example, we have seen this for the case of the symmetrized incomplete gamma functions which occur, for example, when only a finite number of terms in the Riemann A(x) are retained.

The $A^{(1)}(x,k)$ and $A^{(2)}(x,k)$ examples mentioned above both appear to be exceptional because they appear numerically to have a supra-logarithmic growth rate. This suggests the Hurwitz theorem cannot be applied to these cases. The $A^{(1)}(x,k)$ function corresponds to a $\Xi(t)$ which has an infinite number of real zeros, but their distribution is not known exactly. The $A^{(2)}(x,k)$ function has no zeros, and the asymptotics of the corresponding $\Xi(t)$ is not known. The increment equations are therefore difficult to apply, and since in both cases it is known that $\Xi(t)$ has complex zeros off the critical line, we would expect that a future analysis of these equations should result in showing that the growth rate of the minimal beta-sequence cannot be sub-logarithmic.

Historically, the validity of the Riemann Hypothesis has been called into doubt because of the existence of the $A^{(1)}(x,k)$ and $A^{(2)}(x,k)$ cases which have been engineered to show that some Dirichlet series may be similar in many respects to the Riemann series and nevertheless feature complex zeros of their corresponding $\Xi(t)$'s. However, the above analysis of the different reference cases actually show that the Riemann A(x) has much more similarity with simple exponentials than with more complex constructs.

As a final comment on the scope of our analysis of the zeros of the $\Xi(t)$ functions associated with the A(x) family, it is worth emphasizing that the intention of the first few sections concerning the general expansion in symmetrized Pochhammer polynomials, the existence argument for the minimal beta-sequence, and the derivation of the beta and t-increment equations was to provide a clear and mathematically rigorous basis for the subsequent analysis. In the following sections, the intention was more focused on explaining in a simple and intuitive way the fundamental mechanism of our formal analysis and its close agreement with the currently available numerical data, but it would clearly be useful to restate this analysis and its results in a more classical, rigorous form. However, this will not be done in the present paper.

Let us conclude by summarizing a few salient features of our analysis and by speculating on its possible applications in wider contexts than the one discussed in the present paper. The starting point of the analysis can be seen generically as a series or an integral representation for a holomorphic function for which the problem is to elucidate its analytic structure, e.g. the location of complex zeros. If the starting point involves a power function then it is convenient to introduce Pochhammer polynomials into the analysis since their generating function is directly related to this function. At this point, one might imagine other contexts where the applied involution $e \to 1-e$ could be replaced by a different more intricate involution, getting a new type of generating function which in turn would determine the particular possibilities for introducing one or more real dummy parameters, like the $a$ and $b$ used above.

The α and β parameters first appear in terms of an identity, e.g. representing different ways of rearranging terms in an infinite series. However when the series is truncated so as to analyze approximants, then the parameters acquire a more dynamic role which provides for a greater freedom than standard approaches such as multiplier sequences. The question then arises of whether there exists any asymptotic limit of the parameters where the analysis becomes simpler. In the analysis above, we investigated the asymptotic scaling limit of the case $a = b$, however one may note, for example, that some interesting features also appear if β is kept fixed and α becomes asymptotically large. For large $a$ the order 2n polynomial approximants then also happen to enter into a real only root regime where the roots are simple transformations of 2n+1'th roots of unity. There is also here a real root preserving alpha-sequence of the type $a_n = wn$, but in this case the convergence of the approximants turns out be towards the function $\sin(wt)/t$.

It also seems probable that there could be some interesting generalization of the analysis in the direction of using classical orthogonal polynomials in place of Pochhammer polynomials. The orthogonal polynomials all satisfy a similar type of recurrence relations and one observes quite generic features concerning real roots, interlacing sequences, linear combinations and linear transformations. As done for the Pochhammer polynomials, the starting point would be a polynomial $P_k(s)$ with real roots, which is subsequently complexified by a relation similar to

$$P_k(it) = P_k^+(t) + iP_k^-(t) \tag{9.5}$$

in terms of even and odd parts, whereby initially real roots are mapped into other roots on the imaginary axis. This kind of approach might be useful for understanding the even/odd properties of certain analytic functions. The above complexification is explicit in the correspondence between the cosine representation of $\Xi(t)$ and the symmetrized Pochhhammer expansion, as expressed in terms of the polynomial Euler formula (for positive y)

$$\exp(iyt) = \sum_{k=0}^{\infty} e^{-by}(1-e^{-by})^k (P_k^+(t/b) + iP_k^-(t/b)) \tag{9.6}$$

giving rise to the kind of polynomial Fourier transform analysis conducted above. It is possible that this type of polynomial expansion could be useful both numerically and formally not only in different applications involving Fourier transforms, but also Laplace and Mellin transforms.

As a final remark, there may be some reason to believe that the above analysis could be helpful as well for elucidating the various generalized Riemann Hypotheses.

**References**

1. L. V. AHLFORS, *Complex Analysis*, McGraw-Hill, 1966

2. L. BÁEZ-DUARTE, 'A Necessary and Sufficient Condition for the Riemann Hypothesis', arXiv:math-ph/0307215v1, 2003

3. L. BÁEZ-DUARTE, 'On Maślanka's Representation for the Riemann Zeta Function', arXiv:math-ph/0307214v1, 2003

4. S. BELTRAMINELLI AND D. MERLINI, 'Other Representations of the Riemann Zeta Function and an Additional Reformulation of the Riemann Hypothesis', arXiv:math-ph/0707.2406v1, 2007

5. E. BOMBIERI, private communication, 2010

6. J. B. CONREY, 'More than two fifths of the zeros of the Riemann zeta function are on the critical line', *J. reine angew. Math. 399(1989), 1-26.*

7. J. B. CONREY AND A. GHOSH, 'Turán inequalities and zeros of Dirichlet series associated with certain cusp forms', *Trans. Am. Math. Soc. Vol. 342, No.1 (1994)*, and J. B. CONREY, private communication, 2010

8. H. DAVENPORT AND H. HEILBRONN, 'On the zeros of certain Dirichlet series. I, II', *J. London Math. Soc. 11(1936), 1818-85, 307-312*

9. M. W. COFFEY, 'On the Coefficients of the Báez-Duarte criterion for the Riemann Hypothesis and their extensions', arXiv:math-ph/0608050v2, 2006

10. A. M. DIN and L. MONETA, 'Numerical analysis of minimal beta-sequences associated with a family of entire functions', Preprint, July 9, 2011

11. H. M. EDWARDS, *Riemann's Zeta Function*, Dover Publications, 2001

12. S. FENG, 'Zeros of the Riemann zeta function on the critical line', arxiv.org/abs/1003.0059v2, March 2011

13. S. FISK, *Polynomials, Roots, and Interlacing*, Bowdoin College, 2008

14. I. S. GRADSHTEYN and I. M. RYZHIK, *Tables of Integrals, Series and Products*, Academic Press, 1980

15. G. H. HARDY AND J. E. LITTLEWOOD, *Acta Math. 41, 1918, 119-196*

16. A. KUZNETSOV, 'Expansion of the Riemann Xi Function in Meixner-Pollaczek Polynomials', Mc Master University, Canada, 2007

17. B. YA. LEVIN, *Distribution of zeros of entire functions*, English transl., Amer. Math. Soc., 1964

18. J.E. Marsden, *Basic Complex Analysis,* Freeman and Company, 1999

19. K. Maślanka, 'Hypergeometric-like Representation of the Zeta-Function of Riemann', arXiv:math-ph/0105007v1, 2001

20. G. PÓLYA AND G. SZEGÖ, *Problems and Theorems in Analysis Vol.II*, Springer, 1998

21. V. V. P<small>RASOLOV</small>, Polynomials, *Algorithms and Computation in Mathematics*, Springer, 2004

22. M. R<small>IESZ</small>, *Acta Math. 40, 1916, 185-190*

23. A. S<small>ELBERG</small>, 'On the zeros of Riemann's zeta-function', *Skr. Norske Vid. Akad. Oslo.10 (1942),1-59*.

24. E. C. T<small>ITCHMARSH</small>, *The Theory of the Riemann Zeta-Function*, Clarendon Press, 1986

25. L. V<small>EPSTAS</small>, 'A Series Representation of the Riemann Zeta derived from the Gauss-Kuzmin-Wirsing Operator', 2005